\documentclass[draft]{amsart}
\usepackage{amscd}
\usepackage{amssymb}

\newtheorem{thm}[equation]{Theorem}
\newtheorem{cor}[equation]{Corollary}
\newtheorem{prop}[equation]{Proposition}
\newtheorem{lem}[equation]{Lemma}
\theoremstyle{definition}
\newtheorem{dfn}[equation]{Definition}
\newtheorem{rem}[equation]{Remark}
\newtheorem{exa}[equation]{Example}

\newtheorem{prob}[equation]{Problem}
\numberwithin{equation}{section}

\newcommand{\ar}{\rightarrow}
\newcommand{\opn}{\operatorname}
\newcommand{\bdot}{{\textstyle \cdot}}
\newcommand{\ul}{\underline}
\newcommand{\blnk}[1]{\mbox{\hspace{#1}}}
\newcommand{\rmitem}[1]{\item[\text{\textup{(#1)}}]}
\newcommand{\sqbr}[1]{[ \, #1 \, ]}
\newcommand{\mfrak}[1]{\mathfrak{#1}}
\newcommand{\mcal}[1]{\mathcal{#1}}
\newcommand{\msf}[1]{\mathsf{#1}}
\newcommand{\mbf}[1]{\mathbf{#1}}
\newcommand{\mrm}[1]{\mathrm{#1}}
\newcommand{\mbb}[1]{\mathbb{#1}}
\newcommand{\gfrac}[2]{\genfrac{[}{]}{0pt}{}{#1}{#2}}

\title[Smooth Formal Embeddings]{Smooth Formal Embeddings and the
Residue Complex}
\author[Amnon Yekutieli]{Amnon Yekutieli*}
\address{Department of Theoretical Mathematics,
The Weizmann Institute of Science,
Rehovot 76100, ISRAEL}
\date{22 October 1997}
\email{amnon@wisdom.weizmann.ac.il}
\thanks{*Incumbent of the Anna and Maurice Boukstein Career
Development Chair}
\subjclass{Primary: 14B20; Secondary: 14F10, 14B15, 14F20}

\begin{document}

\begin{abstract}
Let $\pi : X \ar S$ be a finite type morphism of noetherian schemes.
A {\em smooth formal embedding} of $X$ (over $S$) is a bijective closed
immersion $X \subset \mfrak{X}$, where $\mfrak{X}$ is a noetherian
formal scheme, formally smooth over $S$.
An example of such an embedding
is the formal completion $\mfrak{X} = Y_{/ X}$ where $X \subset Y$
is an algebraic embedding.
Smooth formal embeddings can be used to calculate algebraic
De Rham (co)homology.

Our main application is an explicit construction of the
Grothendieck residue complex when $S$ is a regular scheme.
By definition the residue complex is the Cousin complex of
$\pi^{!} \mcal{O}_{S}$, as in \cite{RD}.
We start with I-C.\ Huang's theory
of pseudofunctors on modules with $0$-dimensional support,
which provides a graded sheaf
$\bigoplus_{q} \mcal{K}^{q}_{X / S}$.
We then use smooth formal embeddings to obtain the coboundary
operator
$\delta : \mcal{K}^{q}_{X / S} \ar \mcal{K}^{q + 1}_{X / S}$.
We exhibit a canonical isomorphism between the complex
$(\mcal{K}^{\bdot}_{X / S}, \delta)$
and the residue complex of \cite{RD}.
When $\pi$ is equidimensional of dimension $n$ and generically smooth
we show that
$\mrm{H}^{-n} \mcal{K}^{\bdot}_{X / S}$
is canonically isomorphic to to the sheaf of regular differentials
of Kunz-Waldi \cite{KW}.

Another issue we discuss is Grothendieck Duality on a noetherian
formal scheme $\mfrak{X}$.
Our results on duality are used in the construction
of $\mcal{K}^{\bdot}_{X / S}$.
\end{abstract}

\maketitle

% ** section 0 **

\section{Introduction}

It is sometimes the case in algebraic geometry, that in order to define
an object associated to a singular variety $X$, one first embeds $X$
into a nonsingular variety $Y$. One such instance is algebraic
De Rham cohomology
$\mrm{H}^{\bdot}_{\mrm{DR}}(X) =
\mrm{H}^{\bdot}(Y, \widehat{\Omega}^{\bdot})$,
where $\widehat{\Omega}^{\bdot}$ is the completion along $X$ of
the De Rham complex $\Omega_{Y / k}^{\bdot}$ (relative to a base field
$k$ of characteristic $0$; cf.\ \cite{Ha}).
Now $\widehat{\Omega}^{\bdot}$ coincides with the complete De Rham
complex
$\widehat{\Omega}^{\bdot}_{\mfrak{X} / k}$,
where $\mfrak{X}$ is the formal scheme $Y_{/ X}$.
It is therefore reasonable to ask what sort of embedding
$X \subset \mfrak{X}$ into a formal scheme
would give rise to the same cohomology.

The answer we provide in this paper is that any
{\em smooth formal embedding} works. Let us define this notion.
Suppose $S$ is a noetherian base scheme $S$ and
$\pi : X \ar S$ is a finite type morphism.
A smooth formal embedding of $X$ consists of morphisms
$X \ar \mfrak{X} \ar S$, where
$X \ar \mfrak{X}$ is a closed immersion of $X$ into a noetherian
formal scheme $\mfrak{X}$, which is a homeomorphism
of the underlying topological spaces;
and $\mfrak{X} \ar S$ is a {\em formally smooth} morphism.
A smooth formal embedding $X \subset \mfrak{X} = Y_{/ X}$ like
in the previous paragraph is said to be algebraizable. But in general
$X \subset \mfrak{X}$ will not be algebraizable.

Smooth formal embeddings enjoy a few advantages over
algebraic embeddings. First consider an embedding $X \subset \mfrak{X}$
and an \'{e}tale morphism $U \ar X$. Then it is pretty clear
(cf.\ Proposition \ref{prop2.4}) that there is an \'{e}tale
morphism of formal schemes $\mfrak{U} \ar \mfrak{X}$ and a
smooth formal embedding $U \subset \mfrak{U}$, s.t.\
$U \cong \mfrak{U} \times_{\mfrak{X}} X$.
Next suppose $X \subset \mfrak{X}, \mfrak{Y}$ are two smooth formal
embeddings, and we are given either a closed immersion
$\mfrak{X} \ar \mfrak{Y}$ or a formally smooth morphism
$\mfrak{Y} \ar \mfrak{X}$, which restrict to the identity on $X$.
Then locally on $X$,
\begin{equation} \label{eqn0.1}
\mfrak{Y} \cong \mfrak{X} \times \opn{Spf}
\mbb{Z} [\sqbr{t_{1}, \ldots, t_{n}}]
\end{equation}
(Theorem \ref{thm2.2}).

As mentioned above, De Rham cohomology can be calculated by
smooth formal embeddings. Indeed, when $\opn{char} S = 0$,
$\mrm{H}^{q}_{\mrm{DR}}(X / S) =
\mrm{R} \pi_{*} \widehat{\Omega}^{\bdot}_{\mfrak{X} / S}$,
where $X \subset \mfrak{X}$ is any smooth formal embedding
(Corollary \ref{cor2.1}). Moreover, in
\cite{Ye3} it is proved that De Rham homology
$\mrm{H}_{\bdot}^{\mrm{DR}}(X)$
can also be calculated by smooth formal embeddings,
when $S = \opn{Spec} k$, $k$ a field.
According to the preceding paragraph, given an \'{e}tale morphism
$g : U \ar X$ there is a homomorphism
$g^{*} : \mrm{H}_{\bdot}^{\mrm{DR}}(X) \ar
\mrm{H}_{\bdot}^{\mrm{DR}}(U)$,
and we conclude that homology is contravariant w.r.t.\
\'{e}tale morphisms.
See Remark \ref{rem2.4} for an application to $\mcal{D}$-modules
on singular varieties.

The main application of smooth formal embeddings in the present paper
is for an {\em explicit
construction of the Grothendieck residue complex}
$\mcal{K}^{\bdot}_{X / S}$, when $S$ is any regular scheme.
By definition $\mcal{K}^{\bdot}_{X / S}$ is the Cousin complex
$\mrm{E} \pi^{!} \mcal{O}_{S}$, in the notation of \cite{RD}
Sections IV.3 and VII.3.

Recall that Grothendieck Duality, as developed by
Hartshorne in \cite{RD},
is an abstract theory, stated in the language of derived categories.
Even though this abstraction is suitable for many important
applications, often one wants more explicit information.
In particular a significant
amount of work was directed at finding an explicit presentation of
duality in terms of differential forms and residues.
Mostly the focus was on the dualizing sheaf $\omega_{X}$,
in various circumstances. The structure of $\omega_{X}$ as a
coherent $\mcal{O}_{X}$-module and its variance properties are
thoroughly understood by now, thanks to an extended effort
including \cite{KW}, \cite{Li}, \cite{HK1}, \cite{HK2},
\cite{LS1} and \cite{HS}.
Regarding an explicit presentation of the full duality theory of
dualizing complexes, there have been some advances in recent years,
notably in the papers \cite{Ye1}, \cite{SY}, \cite{Hu}, \cite{Hg1}
\cite{Sa} and \cite{Ye3}.
The later papers \cite{Hg2}, \cite{Hg3} and \cite{LS2}
somewhat overlap our present paper in their results, but their methods
are quite distinct; specifically, they do not use formal schemes.

We base our construction of $\mcal{K}^{\bdot}_{X / S}$
on I-C.\ Huang's theory of pseudofunctors on modules with
zero dimensional support (see \cite{Hg1}). Suppose $\phi : A \ar B$
is a residually finitely generated homomorphism between complete
noetherian local rings, and $M$ is a discrete $A$-module
(i.e.\ $\opn{dim} \opn{supp} M = 0$). Then according
to \cite{Hg1} there is a discrete $B$-module $\phi_{\#} M$,
equipped with certain variance properties (cf.\ Theorem \ref{thm6.1}).
In particular
when $\phi$ is residually finite there is a map
$\opn{Tr}_{\phi} : \phi_{\#} M \ar M$. Huang's theory is developed
using only methods of commutative algebra.

Now given a point $x \in X$ with $s := \pi(x) \in S$, consider
the local homomorphism
$\phi : \widehat{\mcal{O}}_{S, s} \ar \widehat{\mcal{O}}_{X, x}$.
Define
$\mcal{K}_{X / S}(x) :=
\phi_{\#} \mrm{H}^{d}_{\mfrak{m}_{s}} \widehat{\mcal{O}}_{S, s}$,
where $d := \opn{dim} \widehat{\mcal{O}}_{S, s}$,
$\mfrak{m}_{s}$ is the maximal ideal and
$\mrm{H}^{d}_{\mfrak{m}_{s}}$ is local cohomology.
Then $\mcal{K}_{X / S}(x)$ is an injective hull of $k(x)$ as
$\mcal{O}_{X, x}$-module.
As a graded $\mcal{O}_{X}$-module we take
$\mcal{K}^{\bdot}_{X / S} := \bigoplus_{x \in X}
\mcal{K}_{X / S}(x)$, with the obvious grading. Then
for any scheme morphism $f : X \ar Y$, we deduce from Huang's theory
a homomorphism of graded sheaves
$\opn{Tr}_{f} : f_{*} \mcal{K}^{\bdot}_{X / S}
\ar \mcal{K}^{\bdot}_{Y / S}$.

The problem is to exhibit a coboundary operator
$\delta : \mcal{K}^{q}_{X / S} \ar \mcal{K}^{q + 1}_{X / S}$,
and to determine that the complex we obtain is indeed isomorphic
to $\mrm{E} \pi^{!} \mcal{O}_{S}$. For this we use smooth formal
embeddings, as explained below.

In Section 5 we discuss Grothendieck Duality on formal schemes,
extending the theory of \cite{RD}.
We propose a definition of dualizing complex $\mcal{R}^{\bdot}$ on
a noetherian formal scheme (Definition \ref{dfn5.1}), and prove its
uniqueness (Theorem \ref{thm5.1}).
It is important to note that the cohomology sheaves
$\mrm{H}^{q} \mcal{R}^{\bdot}$ are discrete quasi-coherent
$\mcal{O}_{\mfrak{X}}$-modules, and in general {\em not coherent}.
We define the Cousin functor
$\mrm{E}$ associated to $\mcal{R}^{\bdot}$,
and show that
$\mrm{E} \mcal{R}^{\bdot} \cong \mcal{R}^{\bdot}$
in the derived category, and
$\mrm{E} \mcal{R}^{\bdot}$ is a residual complex.
On a regular formal scheme $\mfrak{X}$ the (surprising) fact is that
$\mrm{R} \ul{\Gamma}_{\mrm{disc}} \mcal{O}_{\mfrak{X}}$
is a dualizing complex, and not $\mcal{O}_{\mfrak{X}}$
(Theorem \ref{thm5.3}).

Now let $U \subset X$ be an affine open set and suppose
$U \subset \mfrak{U}$ is a smooth formal embedding.
Say $n := \opn{rank} \widehat{\Omega}^{1}_{\mfrak{U} / S}$,
so $\widehat{\Omega}^{n}_{\mfrak{U} / S}$ is a locally free
$\mcal{O}_{\mfrak{U}}$-module of rank $1$,
and
$\mrm{R} \ul{\Gamma}_{\mrm{disc}} \widehat{\Omega}^{n}_{\mfrak{U} / S}[n]$
is a dualizing complex.
Since the Cousin complex is a sum of local cohomology modules,
there is a natural identification of graded
$\mcal{O}_{\mfrak{U}}$-modules
$\mrm{E} \mrm{R} \ul{\Gamma}_{\mrm{disc}}
\widehat{\Omega}^{n}_{\mfrak{U} / S}[n] \cong
\mcal{K}^{\bdot}_{\mfrak{U} / S}$.
This makes
$\mcal{K}^{\bdot}_{\mfrak{U} / S}$ into a complex.
Since
$\mcal{K}^{\bdot}_{U / S} \cong
\mcal{H}om_{\mfrak{U}} \left( \mcal{O}_{U},
\mcal{K}^{\bdot}_{\mfrak{U} / S} \right)$
we come up with an operator $\delta$ on
$\mcal{K}^{\bdot}_{U / S} = \mcal{K}^{\bdot}_{X / S}|_{U}$.

Given another smooth formal embedding $U \subset \mfrak{V}$
we have to compare the complexes $\mcal{K}^{\bdot}_{\mfrak{U} / S}$
and $\mcal{K}^{\bdot}_{\mfrak{V} / S}$. This is rather easy to do using
the following trick.
Choosing a sequence $\ul{a}$ of generators of some defining ideal of
$\mfrak{U}$, and letting $\mbf{K}^{\bdot}_{\infty}(\ul{a})$ be
the associated Koszul complex, we obtain an explicit presentation of
the dualizing complex, namely
\[  \mrm{R} \ul{\Gamma}_{\mrm{disc}}
\widehat{\Omega}^{n}_{\mfrak{U} / S}[n] \cong
\mbf{K}^{\bdot}_{\infty}(\ul{a}) \otimes
\widehat{\Omega}^{n}_{\mfrak{U} / S}[n] \]
(cf.\ Lemma \ref{lem4.3}).
By the structure of smooth formal embeddings
we may assume there is a morphism $f : \mfrak{U} \ar \mfrak{V}$ which is
either formally smooth or a closed immersion. Then choosing relative
coordinates (cf.\ formula \ref{eqn0.1})
and using Koszul complexes we produce a morphism
$\mrm{R} \ul{\Gamma}_{\mrm{disc}} \widehat{\Omega}^{n}_{\mfrak{U} / S}[n]
\ar
\mrm{R} \ul{\Gamma}_{\mrm{disc}} \widehat{\Omega}^{m}_{\mfrak{V} / S}[m]$.
Applying the Cousin functor $\mrm{E}$ we recover
$\opn{Tr}_{f} : \mcal{K}^{\bdot}_{\mfrak{U} / S} \ar
\mcal{K}^{\bdot}_{\mfrak{V} / S}$
as a map of complexes! We conclude that $\delta$ is independent of
$\mfrak{U}$ and hence it glues to a global operator
(Theorem \ref{thm6.2}).

If $f : X \ar Y$ is a finite morphism, then the trace map
$\opn{Tr}_{f} : f_{*} \mcal{K}^{\bdot}_{X / S} \ar
\mcal{K}^{\bdot}_{Y / S}$,
which is provided by Huang's theory, is actually a homomorphism of
complexes (Theorem \ref{thm7.6}).
We show this by employing the same trick as above of going from
Koszul complexes to Cousin complexes, this time inserting a
``Tate residue map'' into the picture.
We use Theorem \ref{thm7.6} to prove directly that if
$\pi : X \ar S$ is equidimensional of dimension $n$ and generically
smooth, then
$\mrm{H}^{-n} \mcal{K}^{\bdot}_{X / S}$
coincides with the sheaf of regular differentials
$\tilde{\omega}^{n}_{X / S}$ of Kunz-Waldi \cite{KW}
(Theorem \ref{thm7.4}).

Finally in Theorem \ref{thm8.10} we exhibit a canonical isomorphism
$\zeta_{X}$
between the complex $\mcal{K}^{\bdot}_{X / S}$ constructed here and
the complex
$\pi^{\triangle} \mcal{O}_{S} = \mrm{E} \pi^{!} \mcal{O}_{S}$
of \cite{RD}.
Given a morphism of schemes $f : X \ar Y$ the isomorphisms
$\zeta_{X}$ and $\zeta_{Y}$ send Huang's trace map
$\opn{Tr}_{f} : f_{*} \mcal{K}^{\bdot}_{X / S} \ar
\mcal{K}^{\bdot}_{Y / S}$
to the trace
$\opn{Tr}^{\mrm{RD}}_{f} : f_{*} \mrm{E} \pi^{!}_{X} \mcal{O}_{S}
\ar \mrm{E} \pi^{!}_{Y} \mcal{O}_{S}$
of \cite{RD} Section VI.4.
In particular it follows that for $f$ proper, $\opn{Tr}_{f}$
is a homomorphism of complexes (Corollary \ref{cor8.1}).

Sections 1 and 3 of the paper contain the necessary supplements
to \cite{EGA}. Perhaps the most noteworthy result there is
Theorem \ref{thm1.10}, which states that formally finite type
morphisms are stable under base change.
This was also proved in \cite{AJL2}.

\medskip \noindent
{\em Acknowledgments.}\
The author wishes to thank L.\ Alonso, I-C.\ Huang, R.\ H\"{u}bl,
A.\ Jerem\'{\i}as, J.\ Lipman and P.\ Sastry for helpful discussions,
some of which took place during a meeting in Oberwolfach in May 1996.

\tableofcontents

% ** section 1 **

\section{Formally Finite Type Morphisms}

In this section we define formally finite type morphisms between
noetherian formal schemes. This mild generalization of the finite
type morphism of \cite{EGA} I \S 10 has the advantage that it
includes the completion morphism
$\mfrak{X} \ar \mfrak{X}_{/ Z}$ (cf.\ Proposition \ref{prop1.12}), and
still is preserved under base change (Theorem \ref{thm1.10}).

We follow the conventions of \cite{EGA} $0_{\mrm{I}}$ \S 7 on {\em adic}
rings. Thus an adic ring is a commutative ring $A$ which is
complete and separated in the $\mfrak{a}$-adic topology, for some ideal
$\mfrak{a} \subset A$.
As for formal schemes, we follow the conventions of \cite{EGA} I
\S 10. Throughout the paper all formal schemes are by default noetherian
(adic) formal schemes.

We write $A \sqbr{\ul{t}} = A \sqbr{t_{1}, \ldots, t_{n}}$
for the polynomial algebra with variables
$t_{1}, \ldots,$ \linebreak
$t_{n}$ over a ring $A$.
The easy lemma below is taken from \cite{AJL2}.

\begin{lem} \label{lem1.11}
Let $A \ar B$ be a continuous homomorphism between noetherian adic rings,
and let $\mfrak{b} \subset B$ be a defining ideal. Then the following
are equivalent:
\begin{enumerate}
\rmitem{i} $A \ar B / \mfrak{b}$ is a finite type homomorphism.
\rmitem{ii} For some homomorphism
$f : A \sqbr{\ul{t}} \ar B$ extending $A \ar B$ one has
$\mfrak{b} = B \cdot f^{-1}(\mfrak{b})$ and
$A \sqbr{\ul{t}} \ar B / \mfrak{b}$ is surjective.
\end{enumerate}
\end{lem}

\begin{proof}
(i) $\Rightarrow$ (ii): Say $b_{1}, \ldots, b_{m}$ generate $\mfrak{b}$
as a $B$-module, and the images of $b_{m+1}, \ldots, b_{n}$ generate
$B / \mfrak{b}$ as an $A$-algebra. Then the homomorphism
$A \sqbr{\ul{t}} \ar B$, $t_{i} \ar b_{i}$ has the required properties.

\noindent (ii) $\Rightarrow$ (i): Trivial.
\end{proof}

\begin{dfn} \label{dfn1.10}
Let $A \ar B$ be a continuous homomorphism between adic noetherian
rings.
We say that $A \ar B$ is of {\em formally finite type} (f.f.t.)
if the equivalent conditions of Lemma \ref{lem1.11} hold.
We shall also say that $B$ is a formally finite type $A$-algebra.
\end{dfn}

\begin{exa} \label{exa1.9}
Let $I \subset A$ be any open ideal, and let
$B := \lim_{\leftarrow i} A / I^{i}$. Then $A \ar B$ is f.f.t.
\end{exa}

Recall that if $A'$ and $B$ are adic $A$-algebras, with defining ideals
$\mfrak{a}'$ and $\mfrak{b}$, the complete tensor product
$A' \widehat{\otimes}_{A} B$ is the completion of
$A' \otimes_{A} B$ w.r.t.\ the topology defined by the image of
$(\mfrak{a}' \otimes_{A} B) \oplus (A' \otimes_{A} \mfrak{b})$.

\begin{prop} \label{prop1.11}
Let $A, A'$ and $B$ be noetherian adic rings, $A \ar B$ a f.f.t.\
homomorphism, and $A \ar A'$ any continuous homomorphism.
Then
$B' := A' \widehat{\otimes}_{A} B$
is a noetherian adic ring, and $A' \ar B'$ is a f.f.t.\ homomorphism.
\end{prop}

\begin{proof}
Choose a homomorphism
$f : A \sqbr{\ul{t}} \ar B$
satisfying condition (ii) of Lemma \ref{lem1.11}. Let
$\mfrak{b} \subset B$ and $\mfrak{a}' \subset A'$ be defining ideals.
Write
$C := A' \otimes_{A} B$ and
$\mfrak{c} := \mfrak{a}' \cdot C + C \cdot \mfrak{b}$,
so
$B' = \lim_{\leftarrow i} C / \mfrak{c}^{i}$.
Consider the homomorphism
$f' : A' \sqbr{\ul{t}} \ar C$, and let
$\mfrak{c}' := {f'}^{-1}(\mfrak{c})$
and
$\widehat{A' \sqbr{\ul{t}}} :=
\lim_{\leftarrow i} A' \sqbr{\ul{t}} / {\mfrak{c}'}^{i}$.
Since $\mfrak{c} = C \cdot \mfrak{c}'$, it follows from
\cite{CA} Section III.2.11 Proposition 14 that
$\widehat{A' \sqbr{\ul{t}}} \ar B'$ is surjective. Hence $B'$ is a
noetherian adic ring with the $\mfrak{b}'$-adic topology, where
$\mfrak{b}' = B' \cdot \mfrak{c}$.
Furthermore
$A' \sqbr{\ul{t}} \ar B' / \mfrak{b}'$
is surjective, and we conclude that $A' \ar B'$ is f.f.t.
\end{proof}

In the next three examples $A$ is an adic ring with defining ideal
$\mfrak{a}$.

\begin{exa} \label{exa1.2}
Recall that for $a \in A$, the complete ring of fractions
$A_{\{a\}}$ is the completion of the localized ring $A_{a}$
w.r.t.\ the $\mfrak{a}_{a}$-adic topology. Then
$A_{\{a\}} \cong A \widehat{\otimes}_{\mbb{Z} \sqbr{t}}
\mbb{Z} \sqbr{t, t^{-1}}$,
which proves that $A \ar A_{\{a\}}$ is f.f.t.
\end{exa}

\begin{exa} \label{exa1.6}
Given indeterminates $t_{1}, \ldots, t_{n}$, the ring of restricted
formal power series
$A\{ \ul{t} \} = A\{ t_{1}, \ldots, t_{n} \}$ is the completion of the
polynomial ring $A \sqbr{ \ul{t} }$ w.r.t.\ the
$(A \sqbr{ \ul{t} } \cdot \mfrak{a})$-adic topology. Hence
$A \{ \ul{t} \} \cong A \widehat{\otimes}_{\mbb{Z}}
\mbb{Z} \sqbr{ \ul{t} }$,
which demonstrates that $A \ar A \{ \ul{t} \}$ is f.f.t.
\end{exa}

\begin{exa} \label{exa1.1}
Consider the adic ring
$A \widehat{\otimes}_{\mbb{Z}} \mbb{Z} [\sqbr{ \ul{t} }]$,
where
$\mbb{Z} [\sqbr{ \ul{t} }] = \mbb{Z} [\sqbr{ t_{1}, \ldots, t_{n} }]$
is the ring of formal power series, with the $(\ul{t})$-adic topology.
Since inverse limits commute, we see that
$A \widehat{\otimes}_{\mbb{Z}} \mbb{Z} [\sqbr{ \ul{t} }]
\cong A [\sqbr{ \ul{t} }]$,
the ring of formal power series over $A$, endowed with the
$(A [\sqbr{ \ul{t} }] \cdot (\mfrak{a} + \ul{t}))$-adic topology.
By Prop.\ \ref{prop1.11},
$A \ar  A [\sqbr{ \ul{t} }]$ is f.f.t.
\end{exa}

Let $A \ar B$ be a f.f.t\ homomorphism between adic rings. Choose
a defining ideal $\mfrak{b} \subset B$, and set
$B_{i} := B / \mfrak{b}^{i+1}$. For $n \geq 0$ define
\[ \widehat{\Omega}^{n}_{B / A} :=
\lim_{\leftarrow i} \Omega^{n}_{B_{i} / A} \cong
\lim_{\leftarrow i} B_{i} \otimes_{B} \Omega^{n}_{B / A} \]
(cf.\ \cite{EGA} $0_{\mrm{IV}}$ 20.7.14).
Let
$\widehat{\Omega}^{\bdot}_{B / A} := \bigoplus_{n \geq 0}
\widehat{\Omega}^{n}_{B / A}$, which is a topological DGA
(differential graded algebra), with
$\widehat{\Omega}^{0}_{B / A} = B$.
This definition is independent of the ideal $\mfrak{b}$.
Since $\Omega^{n}_{B_{i} / A}$ is finite over $B_{i}$ it follows that
$\widehat{\Omega}^{n}_{B / A}$ is finite over $B$.

\begin{rem}
If $A \ar B$ is f.f.t.\ then
$\widehat{\Omega}^{\bdot}_{B/A} \cong \Omega^{\bdot, \mrm{sep}}_{B/A}$,
where $\Omega^{\bdot, \mrm{sep}}_{B/A}$ is the separated algebra
of differentials defined in
\cite{Ye1} \S 1.5 for semi-topo\-logical algebras.
Also $\widehat{\Omega}^{\bdot}_{B/A}$ is the universally finite
differential algebra in the sense of \cite{Ku}.
\end{rem}

\begin{prop} \label{prop1.1}
Let $L \ar A \ar B$ be f.f.t.\ homomorphisms between adic noetherian
rings.
\begin{enumerate}
\item $A \ar B$ is formally smooth relative to $L$ iff  the sequence
\[ 0 \ar B \otimes_{A} \widehat{\Omega}^{1}_{A / L} \xrightarrow{v}
\widehat{\Omega}^{1}_{B / L} \xrightarrow{u}
\widehat{\Omega}^{1}_{B / A} \ar 0 \]
is split exact.
\item  $A \ar B$ is formally \'{e}tale relative to $L$ iff
$B \otimes_{A} \widehat{\Omega}^{1}_{A / L} \ar
\widehat{\Omega}^{1}_{B / L}$
is bijective.
\end{enumerate}
\end{prop}

\begin{proof}
Use the results of \cite{EGA} $0_{\mrm{IV}}$ Section 20.7, together
the fact that these are finite $B$-modules.
\end{proof}

\begin{prop} \label{prop1.3}
Let $f: A \ar B$ be a formally smooth, f.f.t.\  homomorphism
between noetherian adic rings. Then $B$ is flat over $A$ and
$\widehat{\Omega}^{1}_{B/A}$ is a projective finitely generated
$B$-module.
\end{prop}

\begin{proof}
For flatness it suffices to show that if
$\mfrak{n}$ is a maximal ideal of $B$ and
$\mfrak{m} := f^{-1}(\mfrak{n})$,
then
$\widehat{A}_{\mfrak{m}} \ar \widehat{B}_{\mfrak{n}}$
is flat ($\widehat{B}_{\mfrak{n}}$ is the completion of
$B_{\mfrak{n}}$ with the $\mfrak{n}$-adic topology).
Now $\mfrak{n}$ is open,
and hence so is $\mfrak{m}$. Both
$A \ar \widehat{A}_{\mfrak{m}}$ and
$B \ar \widehat{B}_{\mfrak{n}}$ are formally \'{e}tale,
therefore $\widehat{A}_{\mfrak{m}} \ar \widehat{B}_{\mfrak{n}}$
is formally smooth. Because $A \ar B$ is f.f.t.\ it follows that
$A / \mfrak{m} \ar B / \mfrak{n}$
is finite type, and hence finite (and $\mfrak{m}$ is a maximal ideal).
By \cite{EGA} $0_{\mrm{IV}}$ Thm.\ 19.7.1,
$\widehat{B}_{\mfrak{n}}$ is flat over $\widehat{A}_{\mfrak{m}}$.

The second assertion follows from
\cite{EGA} $0_{\mrm{IV}}$ Thm.\ 20.4.9.
\end{proof}

\begin{prop} \label{prop1.4}
Let $f : A \ar B$ be a f.f.t., formally smooth homomorphism of
noetherian adic rings, and let $\mfrak{q} \in \opn{Spf} B$.
Suppose
$\opn{rank} \widehat{\Omega}^{1}_{\widehat{B}_{\mfrak{q}} / A} = n$.
Then:
\begin{enumerate}
\item For some $b \in B - \mfrak{q}$ there is a formally \'{e}tale
homomorphism
$\tilde{f} : A\sqbr{\ul{t}} = A\sqbr{t_{1}, \ldots, t_{n}}
\ar B_{ \{b\} }$ extending $f$.
\item For any
$\mfrak{q}' \in \opn{Spf} B_{ \{b\} }$
let
$\mfrak{r} := \tilde{f}^{-1}(\mfrak{q}')$.
Then
$\widehat{A\sqbr{\ul{t}}}_{\mfrak{r}} \ar \widehat{B}_{\mfrak{q}'}$
is finite \'{e}tale.
\item Let $\mfrak{p} := f^{-1}(\mfrak{q})$.
Assume $\widehat{A}_{\mfrak{p}}$ is regular of dimension $m$,
and
$\opn{tr.deg}_{k(\mfrak{p})} k(\mfrak{q}) = l$. Then
$\widehat{B}_{\mfrak{q}}$ is regular of dimension $n + m - l$.
\end{enumerate}
\end{prop}

\begin{proof}
1.\ By Prop.\ \ref{prop1.3} we can find $b$ s.t.\
$\widehat{\Omega}^{1}_{B_{ \{b\} } / A} \cong B_{ \{b\} } \otimes_{B}
\widehat{\Omega}^{1}_{B / A}$
is free, say with basis $\mrm{d} b_{1}, \ldots, \mrm{d} b_{n}$.
Then we get a homomorphism
$A \sqbr{ \ul{t} } \ar B_{ \{b\} }$, $t_{i} \mapsto b_{i}$.
In order to stay inside the category of adic rings we may replace
$A \sqbr{\ul{t}}$ with its completion $A \{ \ul{t} \}$
(cf.\ Examples \ref{exa1.2} - \ref{exa1.1} for the notation).
According to Proposition \ref{prop1.1} we see that
$A \sqbr{\ul{t}} \ar B_{ \{b\} }$ is formally \'{e}tale relative to $A$.
But since $A \ar B_{ \{b\} }$ is formally smooth, this implies that
$A \sqbr{\ul{t}} \ar B_{ \{b\} }$ is actually (absolutely)
formally \'{e}tale.

\noindent 2.\
Consider the formally \'{e}tale homomorphism
$k(\mfrak{r}) \ar
\widehat{B}_{\mfrak{q}'} / \mfrak{r} \widehat{B}_{\mfrak{q}'}$.
Since $\mfrak{q}'$ is an open prime ideal it follows that
$A \ar B / \mfrak{q}'$ is a finite type homomorphism, so the field
extension $k(\mfrak{r}) \ar k(\mfrak{q}')$ is finitely generated.
By \cite{Hg1} Lemma 3.9 we see that in fact
$\widehat{B}_{\mfrak{q}'} / \mfrak{r} \widehat{B}_{\mfrak{q}'} =
k(\mfrak{q}')$, so
$k(\mfrak{r}) \ar k(\mfrak{q}')$ is finite separable.
Hence
$\widehat{A\sqbr{\ul{t}}}_{\mfrak{r}} \ar \widehat{B}_{\mfrak{q}'}$
is finite \'{e}tale.

\noindent 3.\
Take $\mfrak{q}' := \mfrak{q}$.
Under the assumption the ring
$\widehat{A\sqbr{\ul{t}}}_{\mfrak{r}}$ is regular, and
according to \cite{Ma} \S 14.c Thm.\ 23,
$\opn{dim} \widehat{A\sqbr{\ul{t}}}_{\mfrak{r}} = m + n - l$.
By part 2, $\widehat{B}_{\mfrak{q}}$ is also regular, and
$\opn{dim} \widehat{B}_{\mfrak{q}} =
\opn{dim} \widehat{A\sqbr{\ul{t}}}_{\mfrak{r}}$.
\end{proof}

%%%%

Let us now pass to formal schemes.

Given a noetherian formal scheme $\mfrak{X}$,
choose a defining ideal
$\mcal{I} \subset \mcal{O}_{\mfrak{X}}$, and set
\begin{equation} \label{eqn1.1}
X_{n} := (\mfrak{X}, \mcal{O}_{\mfrak{X}}/\mcal{I}^{n+1}) .
\end{equation}
$X_{n}$ is a noetherian (usual) scheme, and
$\mfrak{X} \cong \lim_{n \ar} X_{n}$
in the category of formal schemes.
One possible choice for $\mcal{I}$ is the largest defining ideal,
in which case one has
$X_{0} = \mfrak{X}_{\mrm{red}}$,
the reduced closed subscheme (cf.\ \cite{EGA} I \S 10.5).

\begin{lem} \label{lem1.12}
Suppose $f : \mfrak{X} \ar \mfrak{Y}$ is a morphism between noetherian
formal schemes. There are defining ideals
$\mcal{I} \subset \mcal{O}_{\mfrak{X}}$ and
$\mcal{J} \subset \mcal{O}_{\mfrak{Y}}$ s.t.\
$f^{-1} \mcal{J} \cdot \mcal{O}_{\mfrak{X}} \subset \mcal{I}$.
Letting $X_{n}$ and $Y_{n}$ be the corresponding schemes
\textup{(}cf.\ \textup{(\ref{eqn1.1})} above\textup{)}, we get morphisms
of schemes
$f_{n} : X_{n} \ar Y_{n}$, with $f = \lim_{n \ar} f_{n}$.
\end{lem}

\begin{proof}
See \cite{EGA} I \S 10.6. For instance, one could take $\mcal{I}$
to be the largest defining ideal and $\mcal{J}$ arbitrary.
\end{proof}

\begin{dfn} \label{dfn1.2}
Let $f: \mfrak{X} \ar \mfrak{Y}$ be a morphism of noetherian (adic)
formal schemes. We say
that $f$ is of {\em formally finite type} (or that $\mfrak{X}$ is a
formally finite type formal scheme over $\mfrak{Y}$) if the morphism
$f_{0} : X_{0} \ar Y_{0}$ in Lemma \ref{lem1.12} is finite type, for some
choice of defining ideals of $\mfrak{X}$ and $\mfrak{Y}$.
\end{dfn}

Observe that if the morphism $f_{0}$ is finite type then so are all
the $f_{n}$, and the definition doesn't depend on the defining ideals
chosen.

\begin{rem} \label{rem1.10}
The definition of f.f.t.\ morphism we gave in an earlier version of the
paper was more cumbersome, though equivalent.
The present Definition \ref{dfn1.2} is taken from \cite{AJL2},
where the name is ``pseudo-finite type morphism'', and I wish to thank
A.\ Jerem\'{\i}as for bringing it to my attention.
\end{rem}

Here are a couple of examples of f.f.t.\ morphisms:

\begin{exa}
A finite type morphism $\mfrak{X} \ar  \mfrak{Y}$ (in the sense of
\cite{EGA} I \S 10.13) is f.f.t.
\end{exa}

\begin{exa} \label{exa1.7}
Let $X$ be a scheme of finite type over a noetherian scheme $S$, and let
$X_{0} \subset X$
be a locally closed subset. Then the completion $\mfrak{X} = X_{/X_{0}}$
(see \cite{EGA} I \S 10.8) is of f.f.t.\
over $S$. Such a formal scheme is called {\em algebraizable}.
\end{exa}

\begin{dfn} \label{dfn1.1}
A f.f.t.\ morphism $f : \mfrak{X} \ar \mfrak{Y}$ is called
{\em formally finite} (resp.\ {\em formally proper}) if the morphism
$f_{0} : X_{0} \ar Y_{0}$ in Lemma \ref{lem1.12} is finite
(resp.\ proper), for some choice of defining ideals.
\end{dfn}

\begin{exa} \label{exa1.8}
If in Example \ref{exa1.7} the subset $X_{0} \subset X$ is closed, then
$\mfrak{X} \ar X$ is formally finite.
If $X_{0} \ar S$ is proper, then $\mfrak{X} \ar S$ is
formally proper.
\end{exa}

\begin{prop} \label{prop1.2}
\begin{enumerate}
\item An immersion $\mfrak{X} \ar \mfrak{Y}$ is f.f.t.
\item If $\mfrak{X} \ar \mfrak{Y}$ and $\mfrak{Y} \ar \mfrak{Z}$ are
f.f.t., then so is $\mfrak{X} \ar \mfrak{Z}$.
\item Let $\mfrak{U} = \opn{Spf} B$ and
$\mfrak{V} = \opn{Spf} A$. Then a morphism
$\mfrak{U} \ar \mfrak{V}$ is f.f.t.\ iff the ring homomorphism
$A \ar B$ is f.f.t.
\end{enumerate}
\end{prop}

\begin{proof}
Consider morphisms of schemes $X_{0} \ar Y_{0}$ etc.\
as in Lemma \ref{lem1.12}.
For part 3 use condition (i) of Lemma \ref{lem1.11}.
\end{proof}

\begin{prop} \label{prop1.12}
Let $\mfrak{X}$ be a noetherian formal scheme and $Z \subset \mfrak{X}$
a locally closed subset. Then there is a noetherian formal scheme
$\mfrak{X}_{/Z}$, with underlying topological space $Z$, and the
natural morphism $\mfrak{X}_{/Z} \ar \mfrak{X}$
is f.f.t.
\end{prop}

\begin{proof}
Pick an open subset $\mfrak{U} \subset \mfrak{X}$ s.t.\
$Z \subset \mfrak{U}$ is closed, and choose a defining ideal
$\mcal{I}$ of $Z$. Let
$\mcal{O}_{\mfrak{Z}} := \lim_{\leftarrow i}
\mcal{O}_{\mfrak{U}} / \mcal{I}^{i}$.
According to \cite{EGA} I Section 10.6,
$\mfrak{X}_{/Z} := (Z, \mcal{O}_{\mfrak{Z}})$
is a noetherian formal scheme. Clearly
$\mfrak{X}_{/Z} \ar \mfrak{X}$ is f.f.t.
\end{proof}

In \cite{EGA} I \S 10.3 it is shown that finite type morphisms
between noetherian formal schemes are preserved by base change.
This is true also for f.f.t.\ morphisms:

\begin{thm} \label{thm1.10}
Suppose $\mfrak{X}$, $\mfrak{Y}$ and $\mfrak{Y}'$ are noetherian
formal schemes,
$\mfrak{X} \ar \mfrak{Y}$ is a f.f.t.\ morphism, and
$\mfrak{Y}' \ar \mfrak{Y}$ is an arbitrary morphism. Then
$\mfrak{X}' := \mfrak{X} \times_{\mfrak{Y}} \mfrak{Y}'$
is also noetherian, and the morphism $\mfrak{X}' \ar \mfrak{Y}'$
is f.f.t.
\end{thm}

\begin{proof}
First note that the formal scheme
$\mfrak{X}' = \mfrak{X} \times_{\mfrak{Y}} \mfrak{Y}'$
exists (\cite{EGA} I \S 10.7).
For any affine open sets
$\mfrak{U} = \opn{Spf} B \subset \mfrak{X}$,
$\mfrak{V}' = \opn{Spf} A' \subset \mfrak{Y}'$ and
$\mfrak{V} = \opn{Spf} A \subset \mfrak{Y}$ such that
$\mfrak{U} \ar \mfrak{V}$ and $\mfrak{V}' \ar \mfrak{V}$,
one has
$\mfrak{U}' = \mfrak{U} \times_{\mfrak{V}} \mfrak{V}' =
\opn{Spf} B \widehat{\otimes}_{A} A'$,
and $\mfrak{U}' \subset \mfrak{X}'$ is open.
By Propositions \ref{prop1.11} and \ref{prop1.2}, $\mfrak{U}'$ is a
noetherian formal scheme, and
$\mfrak{U}' \ar \mfrak{V}'$ is f.f.t. But finitely many such
$\mfrak{U}'$ cover $\mfrak{X}'$.
\end{proof}

\begin{cor} \label{cor1.10}
If $\mfrak{X}_{1}$, $\mfrak{X}_{2}$ and $\mfrak{Y}$ are
noetherian and $\mfrak{X}_{i} \ar \mfrak{Y}$ are f.f.t.\ morphisms,
then
$\mfrak{X}_{3} := \mfrak{X}_{1} \times_{\mfrak{Y}} \mfrak{X}_{2}$
is also noetherian, and
$\mfrak{X}_{3} \ar \mfrak{Y}$ is f.f.t.
\end{cor}

\begin{rem} \label{rem1.4}
I do not know an example of a f.f.t.\ formal scheme $\mfrak{X}$ over a
scheme $S$ which is not locally algebraizable. (Locally algebraizable
means there is an open
covering $\mfrak{X} = \bigcup \mfrak{U}_{i}$, with $\mfrak{U}_{i} \ar S$
algebraizable, in the sense of Example \ref{exa1.7}.)
\end{rem}

\begin{dfn} \label{dfn1.3}
A morphism of formal schemes $\mfrak{X} \ar \mfrak{Y}$ is said to be {\em
formally smooth} (resp.\ {\em formally \'{e}tale}) if, given
a (usual) affine scheme $Z$ and a closed subscheme $Z_{0} \subset Z$
defined by a nilpotent ideal, the map
$\opn{Hom}_{\mfrak{Y}}(Z, \mfrak{X}) \ar
\opn{Hom}_{\mfrak{Y}}(Z_{0}, \mfrak{X})$
is surjective (resp.\ bijective).
\end{dfn}

This is the definition of formal smoothness used in \cite{EGA} IV
Section 17.1. We shall also require the next notion.

\begin{dfn} \label{dfn1.4}
A morphism $g: \mfrak{X} \ar \mfrak{Y}$ between noetherian formal schemes
is called {\em \'{e}tale} if it is of finite
type (see \cite{EGA} I \S 10.13) and formally \'{e}tale.
\end{dfn}

Note that if $\mfrak{Y}$ is a usual scheme, then so is
$\mfrak{X}$, and $g$ is an \'{e}tale morphism of schemes.
According to \cite{EGA} I Prop.\ 10.13.5 and by the obvious properties
of formally \'{e}tale morphisms, if $\mfrak{U} \ar \mfrak{X}$ and
$\mfrak{V} \ar \mfrak{X}$ are \'{e}tale, then so is
$\mfrak{U} \times_{\mfrak{X}} \mfrak{V} \ar \mfrak{X}$.
Hence for fixed $\mfrak{X}$, the
category of all \'{e}tale morphisms $\mfrak{U} \ar \mfrak{X}$ forms a site
(cf.\ \cite{Mi} Ch.\ II \S 1). We call this site the small \'{e}tale site
on $\mfrak{X}$, and denote it by $\mfrak{X}_{\mrm{et}}$.

% ** section 2 **

\section{Smooth Formal Embeddings and De Rham Cohomology}

Fix a noetherian base scheme $S$ and a finite type $S$-scheme $X$.

\begin{dfn} \label{dfn2.1}
A {\em smooth formal embedding} (s.f.e.) of $X$ (over $S$)
is the following data:
\begin{enumerate}
\rmitem{i} A noetherian formal scheme $\mfrak{X}$.
\rmitem{ii} A formally finite type, formally smooth morphism
$\mfrak{X} \ar S$.
\rmitem{iii} An $S$-morphism $X \ar \mfrak{X}$, which is a
closed immersion and a homeomorphism between the underlying
topological spaces.
\end{enumerate}
We shall refer to this by writing ``$X \subset \mfrak{X}$ is a s.f.e.''
\end{dfn}

\begin{exa} \label{exa2.1}
Suppose $Y$ is a smooth $S$-scheme, $X \subset Y$ a locally closed subset,
and $\mfrak{X} = Y_{/X}$ the completion. Then $X \subset \mfrak{X}$ is a
smooth formal embedding. Such an embedding is called an
{\em algebraizable embedding} (cf.\ Remark \ref{rem1.4}).
\end{exa}

The smooth formal embeddings of $X$ form a category, in which a morphism
of embeddings is an $S$-morphism of formal schemes
$f : \mfrak{X} \ar \mfrak{Y}$
inducing the identity on $X$.
Note that any morphism of embeddings $f: \mfrak{X} \ar \mfrak{Y}$
is affine (cf.\ \cite{EGA} I Prop.\ 10.6.12), and the functor
$f_{*} : \mathsf{Mod}(\mfrak{X}) \ar \mathsf{Mod}(\mfrak{Y})$
is exact. Let
$\mfrak{X}$ and $\mfrak{Y}$ be two smooth formal embeddings of $X$.
Consider the formal scheme
$\mfrak{X} \times_{S} \mfrak{Y}$. Then the diagonal
$\Delta : X \ar \mfrak{X} \times_{S} \mfrak{Y}$
is an immersion (we do not assume our formal schemes are separated!).

\begin{prop} \label{prop2.1}
The completion
$(\mfrak{X} \times_{S} \mfrak{Y})_{/ X}$
of $\mfrak{X} \times_{S} \mfrak{Y}$ along $\Delta(X)$ is a smooth formal
embedding of $X$, and moreover it is a product of $\mfrak{X}$ and
$\mfrak{Y}$ in the category of smooth formal embeddings.
\end{prop}

\begin{proof}
By Theorem \ref{thm1.10} and Proposition \ref{prop1.12} it follows that
$(\mfrak{X} \times_{S} \mfrak{Y})_{/ X}$ is
formally finite type over $S$, so in particular it is noetherian.
Clearly
$(\mfrak{X} \times_{S} \mfrak{Y})_{/ X} \ar S$
is formally smooth.
\end{proof}

The benefit of using formal rather than algebraic embeddings is in:

\begin{prop} \label{prop2.4}
Let $X \subset \mfrak{X}$ be a smooth formal embedding \textup{(}over
$S$\textup{)}
and $g : U \ar X$ an \'{e}tale morphism. Then there exists a noetherian
formal scheme $\mfrak{U}$ and an \'{e}tale morphism
$\widehat{g} : \mfrak{U} \ar \mfrak{X}$ s.t.\
$U \cong \mfrak{U} \times_{\mfrak{X}} X$.
$\widehat{g} : \mfrak{U} \ar \mfrak{X}$
is unique \textup{(}up to a unique isomorphism\textup{)}, and moreover
$U \ar \mfrak{U}$ is a smooth formal embedding.
\end{prop}

\begin{proof}
This is essentially the ``topological invariance of \'{e}tale morphisms'',
cf.\ \cite{EGA} IV \S 18.1 (or \cite{Mi} Ch.\ I Thm.\ 3.23).
Let
$\mcal{I} := \opn{Ker}(\mcal{O}_{\mfrak{X}} \ar \mcal{O}_{X})$ and
$X_{i} := (\mfrak{X}, \mcal{O}_{\mfrak{X}} / \mcal{I}^{i+1})$;
so $X = X_{0}$. For every $i$ there is a unique \'{e}tale morphism
$g_{i} : U_{i} \ar X_{i}$ s.t.\
$U \cong U_{i} \times_{X_{i}} X$. Identifying the underlying topological
spaces of $U_{i}$ and $U$ we get an inverse system of sheaves
$\{ \mcal{O}_{U_{i}} \}$ on $U$. Setting
$\mcal{O}_{\mfrak{U}} := \lim_{\leftarrow i} \mcal{O}_{U_{i}}$ we get a
noetherian formal scheme $\mfrak{U}$ having the proclaimed properties
(cf.\ \cite{EGA} I \S 10.6).
\end{proof}

Thus we can consider $\mfrak{X}_{\mrm{et}}$ as a ``smooth formal
embedding'' of $X_{\mrm{et}}$.
If $\mcal{M}$ is a sheaf on $X_{\mrm{et}}$ and $U \ar X$ is an \'{e}tale
morphism, we denote by $\mcal{M}|_{U}$ the restriction of $\mcal{M}$ to
$U_{\mrm{Zar}}$.

\begin{cor} \label{cor2.3}
Let $X \subset \mfrak{X}$ be a smooth formal embedding over $S$.
Then there is a sheaf of DGAs
$\widehat{\Omega}^{\bdot}_{\mfrak{X}_{\mrm{et}} / S}$
on $X_{\mrm{et}}$, with the property that for each
$g: U \ar X$ in $X_{\mrm{et}}$ and corresponding
$\widehat{g}: \mfrak{U} \ar \mfrak{X}$ in $\mfrak{X}_{\mrm{et}}$, one has
$\widehat{\Omega}^{\bdot}_{\mfrak{X}_{\mrm{et}} / S} |_{U} \cong
\widehat{\Omega}^{\bdot}_{\mfrak{U} / S}$.
\end{cor}

\begin{proof}
By Prop.\ \ref{prop1.1},
$\widehat{\Omega}^{p}_{\mfrak{U} / S} \cong
\widehat{g}^{*} \widehat{\Omega}^{p}_{\mfrak{X} / S}$.
Now $\widehat{\Omega}^{p}_{\mfrak{X} / S}$ is coherent, so we can use
\cite{Mi} Ch.\ II Cor.\ 1.6 (which applies to our \'{e}tale site
$\mfrak{X}_{\mrm{et}}$).
\end{proof}

For smooth formal embeddings, closed immersions and smooth morphisms
are locally trivial, in the following sense. Recall that for an
adic algebra $A$, the ring of formal power series
$A [\sqbr{ \ul{t} }] = A [\sqbr{ t_{1}, \ldots, t_{n} }]$
is adic (cf.\ Example \ref{exa1.1}).

\begin{thm} \label{thm2.2}
Let $f : \mfrak{X} \ar \mfrak{Y}$ be a morphism of smooth formal
embeddings of $X$ over $S$. Assume $f$ is a closed immersion
\textup{(}resp.\ formally smooth\textup{)}.
Then, given a point $x \in X$, there are affine open sets $U \subset X$
and $W \subset S$, with $x \in U$ and $U \ar W$, satisfying condition
\textup{($*$)} below.
\begin{enumerate}
\item[\textup{($*$)}]
Let $W = \opn{Spec} L$, and let
$\opn{Spf} A \subset \mfrak{Y}$ and
$\opn{Spf} B \subset \mfrak{X}$
be the affine formal schemes supported on $U$. Then
there is an isomorphism of topological $L$-algebras
$A \cong B [\sqbr{\ul{t}}]$
\textup{(}resp.\ $B \cong A [\sqbr{\ul{t}}]$\textup{)}
such that $f^{*} : A \ar B$ is projection modulo $(\ul{t})$
\textup{(}resp.\ the inclusion\textup{)}.
\end{enumerate}
\end{thm}

\begin{proof}
1.\ Assume $f$ is a closed immersion.
According to \cite{EGA} $0_{\mrm{IV}}$ Thm.\ 19.5.3 and
Cor.\ 20.7.9, by choosing $U = \opn{Spec} C$ small enough, and setting
$I := \opn{Ker}(f^{*} : A \ar B)$, we obtain
an exact sequence
\[ 0 \ar I / I^{2} \ar B \otimes_{A} \widehat{\Omega}^{1}_{A/L} \ar
\widehat{\Omega}^{1}_{B/L} \ar 0 \]
of free $B$-modules. Choose
$a_{1}, \ldots, a_{n}, b_{1}, \ldots, b_{m} \in A$
s.t.\ $\{ a_{i} \}$ is a $B$-basis of $I / I^{2}$, and
$\{ \mrm{d} b_{i} \}$
is a $B$-basis of $\widehat{\Omega}^{1}_{B/L}$.

By the proof of Prop.\ \ref{prop1.4} the homomorphisms
$L \sqbr{\ul{s}} \ar B$,
$L \sqbr{\ul{s}, \ul{t}} \ar A$ and
$L \sqbr{\ul{s}, \ul{t}} \ar B [\sqbr{\ul{t}}]$,
sending
$s_{i} \mapsto b_{i}$ and $t_{i} \mapsto a_{i}$,
are all formally \'{e}tale.
Take $\mfrak{a} := \opn{Ker}(A \ar C)$, which is a defining ideal of
$A$, containing $A \cdot (\ul{t}) = I$.
Let $\mfrak{b} := \mfrak{a} \cdot B$, which is a
defining ideal of $B$. Hence the ideal
$\mfrak{c} = B [\sqbr{\ul{t}}] \cdot (\mfrak{b}, \ul{t})$
is a defining ideal of $B [\sqbr{\ul{t}}]$. By formal \'{e}taleness
of
$L \sqbr{\ul{s}, \ul{t}} \ar A$ and
$L \sqbr{\ul{s}, \ul{t}} \ar B [\sqbr{\ul{t}}]$,
the isomorphism
$A / \mfrak{a} \cong B [\sqbr{\ul{t}}] / \mfrak{c} \cong C$
lifts uniquely to an isomorphism
$A \cong B [\sqbr{\ul{t}}]$.

\noindent 2.\
Now assume $f$ is formally smooth. Let
$\mfrak{b} := \opn{Ker}(B \ar C)$,
which is a defining ideal of $B$.
Since $A \ar B / \mfrak{b}$ is surjective it follows that
$(B / \mfrak{b}) \otimes_{B} \widehat{\Omega}^{1}_{B/A}$
is generated by $\mrm{d}(\mfrak{b})$. By Nakayama's Lemma we see that
$\widehat{\Omega}^{1}_{B/A} = B \cdot \mrm{d}(\mfrak{b})$.
After shrinking $U$ sufficiently we get
$\widehat{\Omega}^{1}_{B/A} =
\bigoplus_{i = 1}^{n} B \cdot \mrm{d} b_{i}$
with $b_{i} \in \mfrak{b}$, and the homomorphism
$A [\sqbr{\ul{t}}] \ar B$, $t_{i} \mapsto b_{i}$,
is formally \'{e}tale. Continuing like in part 1 of the proof we
conclude that this is actually an isomorphism.
\end{proof}

\begin{thm} \label{thm2.1}
Suppose $S$ is a noetherian scheme of characteristic $0$, and $X$ is a
finite type $S$-scheme.
Let $f : \mfrak{X} \ar \mfrak{Y}$ be a morphism of smooth formal
embeddings of $X$. Then the DGA homomorphism
$f^{*} : \widehat{\Omega}^{\bdot}_{\mfrak{Y} / S} \ar
\widehat{\Omega}^{\bdot}_{\mfrak{X} / S}$
is a quasi-isomorphism. Moreover, if $g : \mfrak{X} \ar \mfrak{Y}$ is any
other morphism, then
$\mrm{H}(f^{*}) = \mrm{H}(g^{*})$.
\end{thm}

\begin{proof}
The assertions of the theorem are both local, and they will be proved in
three steps.

\noindent Step 1.\
Assume $f$ is a closed immersion. By Thm.\ \ref{thm2.2}
it suffices to check the case
$f  : \opn{Spf} B = \mfrak{U} \ar \opn{Spf} A = \mfrak{V}$
with $A \cong B [\sqbr{\ul{t}}]$ as topological $L$-algebras.
We must show that
$\widehat{\Omega}^{\bdot}_{A/L} \ar \widehat{\Omega}^{\bdot}_{B/L}$
is a quasi-isomorphism. But since $\mbb{Q} \subset L$,
this is the well known Poincar\'{e} Lemma for
formal power series (cf.\ \cite{Ha} Ch.\ II Prop.\ 1.1, or \cite{Ye3}
Lemma 7.5).

\noindent Step 2.\
Suppose $f_{1}, f_{2}: \mfrak{X} \ar \mfrak{Y}$ are two morphisms.
We wish to show that $\mrm{H}(f_{1}^{*}) = \mrm{H}(f_{2}^{*})$.
First consider
\[ \mfrak{Y} \xrightarrow{\mrm{diag}} (\mfrak{Y} \times_{k} \mfrak{Y})_{/X}
\xrightarrow{p_{i}} \mfrak{Y} . \]
Since the diagonal immersion is closed, we can apply the result of the
previous paragraph to it. We conclude that
$\mrm{H}(p_{1}^{*}) = \mrm{H}(p_{2}^{*})$,
and that these are isomorphisms.
But looking at
\[ \mfrak{X} \xrightarrow{\mrm{diag}} (\mfrak{X} \times_{k}
\mfrak{X})_{/X}
\xrightarrow{f_{1} \times f_{2}} (\mfrak{Y} \times_{k} \mfrak{Y})_{/X}
\xrightarrow{p_{i}} \mfrak{Y} \]
we see that our claim is proved.

\noindent Step 3.\
Consider an arbitrary morphism
$f: \mfrak{X} \ar \mfrak{Y}$. Take any affine open
set $U \subset X$, with corresponding affine formal schemes
$\opn{Spf} B = \mfrak{U} \subset \mfrak{X}$ and
$\opn{Spf} A = \mfrak{V} \subset \mfrak{Y}$.
The definition of formal smoothness implies
there is some morphism of embeddings $g : \mfrak{V} \ar \mfrak{U}$. This
morphism will not necessarily be an inverse of $f|_{\mfrak{U}}$, but
nonetheless, according to Step 2, $\mrm{H}(g^{*})$ and
$\mrm{H}(f|_{\mfrak{U}}^{*})$ will be isomorphisms between
$\mrm{H} \widehat{\Omega}^{\bdot}_{\mfrak{U} / S}$
and
$\mrm{H} \widehat{\Omega}^{\bdot}_{\mfrak{V} / S}$, inverse to each
other.
\end{proof}

In \cite{Ha} the relative De Rham cohomology
$\mrm{H}^{\bdot}_{\mrm{DR}} (X / S)$
was defined. In the situation of Example \ref{exa2.1},
where $X \subset Y$ is a smooth algebraic embedding of $S$-schemes,
$\mfrak{X} = Y_{/X}$ and $\pi : \mfrak{X} \ar S$ is the structural
morphism, the definition is
$\mrm{H}^{\bdot}_{\mrm{DR}} (X / S) =
\mrm{H}^{\bdot} \mrm{R} \pi_{*}
\widehat{\Omega}^{\bdot}_{\mfrak{X} / S}$.
Even if $X$ is not globally embeddable,
$\mrm{H}^{\bdot}_{\mrm{DR}} (X / S)$ can still be defined, by taking
a system of local embeddings
$\{ U_{i} \subset V_{i} \}$, $X = \bigcup U_{i}$,
and putting together a ``\v{C}ech-De Rham'' complex (cf.\ \cite{Ha} pp.\
28-29; it seems one should also assume $X$ separated and the
$U_{i}$ are affine).

\begin{cor} \label{cor2.1}
Suppose $S$ has characteristic $0$.
Let $X \subset \mfrak{X}$ be any smooth formal embedding \textup{(}not
necessarily algebraizable\textup{)}. Then
$\mrm{H}^{\bdot}_{\mrm{DR}} (X / S) = \mrm{H}^{\bdot} \mrm{R} \pi_{*}
\widehat{\Omega}^{\bdot}_{\mfrak{X} / S}$
as graded $\mcal{O}_{S}$-algebras.
\end{cor}

\begin{proof}
Assume for simplicity that a global smooth algebraic embedding exists.
The general case, involving a system of embeddings,
only requires more bookkeeping.
Say $X \subset Y$ is the given algebraic embedding, and let
$\mfrak{Y} := Y_{/X}$.
Now the two formal embeddings $\mfrak{X}$ and $\mfrak{Y}$
are comparable: their product $(\mfrak{X} \times_{S} \mfrak{Y})_{/X}$
maps to both. By the theorem we get quasi-isomorphic DGAs on $X$.
\end{proof}

\begin{rem} \label{rem2.1}
From Corollaries \ref{cor2.3} and \ref{cor2.1} we see that there is a
sheaf of DGAs
$\widehat{\Omega}^{\bdot}_{\mfrak{X}_{\mrm{et}} / S}$ on
$X_{\mrm{et}}$, with the property that for any $U \ar X$ \'{e}tale,
$\mrm{H}^{\bdot}_{\mrm{DR}}(U/S) =$ \linebreak
$\mrm{H}^{\bdot} \Gamma(U,
\widehat{\Omega}^{\bdot}_{\mfrak{X}_{\mrm{et}} / S})$.
As will be shown in \cite{Ye4}, the DGA
$\widehat{\Omega}^{\bdot}_{\mfrak{X} / S}$ has an adelic resolution
$\mcal{A}^{\bdot}_{\mfrak{X} / S}$,
where
$\mcal{A}^{p,q}_{\mfrak{X} / S} = \ul{\mbb{A}}^{q}_{\mrm{red}}(
\widehat{\Omega}^{p}_{\mfrak{X} / S})$, Beilinson's sheaf of adeles.
The adeles calculate cohomology:
$\mrm{H}^{\bdot}_{\mrm{DR}}(X/S) =
\mrm{H}^{\bdot} \Gamma(X, \mcal{A}^{\bdot}_{\mfrak{X} / S})$.
Furthermore the adeles extend to an \'{e}tale sheaf
$\mcal{A}^{\bdot}_{\mfrak{X}_{\mrm{et}} / S}$.
\end{rem}

\begin{rem} \label{rem2.3}
Suppose $S = \opn{Spec} k$, a field of characteristic $0$.
In \cite{Ye3} a complex
$\mcal{F}^{\bdot}_{\mfrak{X}}$, called the De Rham-residue complex, is
defined. One has
$\mrm{H}^{i}(X, \mcal{F}^{\bdot}_{\mfrak{X}}) =
\mrm{H}_{-i}^{\mrm{DR}}(X)$, the De Rham homology.
Moreover there is a sheaf
$\mcal{F}^{\bdot}_{\mfrak{X}_{\mrm{et}}}$
on $X_{\mrm{et}}$, which directly implies that the De Rham homology
is contravariant for \'{e}tale morphisms.
Furthermore $\mcal{F}^{\bdot}_{\mfrak{X}}$
is naturally a DG $\mcal{A}^{\bdot}_{\mfrak{X}}$-module,
\end{rem}

\begin{rem} \label{rem2.4}
Smooth formal embeddings can be also used to define the category
of $\mcal{D}$-modules on a singular scheme $X$ (in characteristic $0$).
Say $X \subset \mfrak{X}$ is such an embedding. Then a formal
version of Kashiwara's Theorem (cf.\ \cite{Bo} Theorem VI.7.11)
implies that
$\msf{Mod}_{\mrm{disc}}(\mcal{D}_{\mfrak{X}})$,
the category of discrete modules over the ring of differential operators
$\mcal{D}_{\mfrak{X}}$ is, as an abelian category, independent of
$\mfrak{X}$.
\end{rem}

% ** section 3 **

\section{Quasi-Coherent Sheaves on Formal Schemes}

Let $\mfrak{X}$ be a noetherian (adic) formal scheme.
By definition, a quasi-coherent sheaf on $\mfrak{X}$ is an
$\mcal{O}_{\mfrak{X}}$-module $\mcal{M}$, such that on sufficiently small
open sets $\mfrak{U} \subset \mfrak{X}$  there are exact sequences
$\mcal{O}_{\mfrak{U}}^{(J)} \ar \mcal{O}_{\mfrak{U}}^{(I)} \ar
\mcal{M}|_{\mfrak{U}} \ar 0$,
for some indexing sets $I,J$ (cf.\ \cite{EGA} $0_{\mrm{I}}$ \S 5.1).
We shall denote by
$\mathsf{Mod}(\mfrak{X})$ (resp.\ $\mathsf{Coh}(\mfrak{X})$, resp.\
$\mathsf{QCo}(\mfrak{X})$)
the category of $\mcal{O}_{\mfrak{X}}$-modules (resp.\ the full subcategory
of coherent, resp.\ quasi-coherent, modules).
It seems that the only important quasi-coherent sheaves are the coherent
and the discrete ones (Def.\ \ref{dfn3.1}). Nevertheless we shall
consider all quasi-coherent sheaves, at the price of a little extra
effort.

\begin{rem}
There is some overlap between results in this section and
\cite{AJL2}.
\end{rem}

Let $A$ be a noetherian adic ring, and let
$\mfrak{U} := \opn{Spf} A$ be the affine formal scheme.
Then there is an exact functor $M \mapsto M^{\triangle}$ from the category
$\mathsf{Mod}_{\mrm{f}}(A)$ of
finitely generated $A$-modules to $\mathsf{Mod}(\mfrak{U})$.
It is an equivalence between $\mathsf{Mod}_{\mrm{f}}(A)$ and
$\mathsf{Coh}(\mfrak{U})$ (see \cite{EGA} I \S 10.10).

\begin{prop} \label{prop3.2}
The functor $M \mapsto M^{\triangle}$ extends uniquely to a functor
$\mathsf{Mod}(A) \ar \mathsf{Mod}(\mfrak{U})$,
which is exact and commutes with direct limits.
The $\mcal{O}_{\mfrak{U}}$-module $M^{\triangle}$ is
quasi-coherent. For any $\mcal{O}_{\mfrak{U}}$-module $\mcal{M}$ the
following are equivalent:
\begin{enumerate}
\rmitem{i} $\mcal{M} \cong M^{\triangle}$ for some $A$-module $M$.
\rmitem{ii} $\mcal{M} \cong \lim_{\alpha \ar} \mcal{M}_{\alpha}$
for some directed system $\{ \mcal{M}_{\alpha} \}$ of coherent
$\mcal{O}_{\mfrak{U}}$-modules.
\rmitem{iii} For every affine open set
$\mfrak{V} = \opn{Spf} B \subset \mfrak{U}$, one has
$\Gamma(\mfrak{V}, \mcal{M}) \cong
B \otimes_{A} \Gamma(\mfrak{U}, \mcal{M})$.
\end{enumerate}
\end{prop}

\begin{proof}
Take any $A$-module $M$ and write it as
$M = \lim_{\alpha \ar} M_{\alpha}$ with finitely generated modules
$M_{\alpha}$. Define a presheaf $M^{\triangle}$ on $\mfrak{U}$ by
$\Gamma(\mfrak{V}, M^{\triangle}) := \lim_{\alpha \ar}
\Gamma(\mfrak{V}, M^{\triangle}_{\alpha})$, for
$\mfrak{V} \subset \mfrak{U}$
open. Since $\mfrak{U}$ is a noetherian topological space it follows that
$M^{\triangle}$ is actually a sheaf. By construction
$M \mapsto M^{\triangle}$ commutes with direct limits. Since the
functor is exact on $\mathsf{Mod}_{\mrm{f}}(A)$, it's also exact on
$\mathsf{Mod}(A)$.

The implication (i) $\Rightarrow$ (ii) is because $M^{\triangle}_{\alpha}$
is coherent.
(ii) $\Rightarrow$ (iii): for such $B$ one has
$\Gamma(\mfrak{V}, \mcal{M}_{\alpha}) \cong
B \otimes_{A} \Gamma(\mfrak{U}, \mcal{M}_{\alpha})$;
now apply $\lim_{\alpha \ar}$.
(iii) $\Rightarrow$ (i): set
$M := \Gamma(\mfrak{U}, \mcal{M})$. Then for every affine $\mfrak{V}$ we
have
$\Gamma(\mfrak{V}, \mcal{M}) = B \otimes_{A} M =
\Gamma(\mfrak{V}, M^{\triangle})$,
so $\mcal{M} = M^{\triangle}$.

Finally the module $M$ has a presentation
$A^{(I)} \ar A^{(J)} \ar M \ar 0$. By exactness we get a presentation for
$M^{\triangle}$.
\end{proof}

It will be convenient to write $\mcal{O}_{\mfrak{U}} \otimes_{A} M$
instead of $M^{\triangle}$.

\begin{rem}
I do not know whether Serre's Theorem holds,
namely whe\-ther {\em every} quasi-coherent $\mcal{O}_{\mfrak{U}}$-module
$\mcal{M}$ is of the form
$\mcal{M} \cong \mcal{O}_{\mfrak{U}} \otimes_{A} M$.
Thus it may be that $\mathsf{QCo}(\mfrak{U})$ is not closed under direct
limits in $\mathsf{Mod}(\mfrak{U})$ (cf.\ Lemma \ref{lem4.1}).
\end{rem}

\begin{cor} \label{cor3.1}
Let $\mcal{M}$ be a quasi-coherent $\mcal{O}_{\mfrak{X}}$-module and
$x \in \mfrak{X}$ a point. Then there is an open neighborhood
$\mfrak{U} = \opn{Spf} A$ of $x$ s.t.\
$\mcal{M}|_{\mfrak{U}} \cong \mcal{O}_{\mfrak{U}} \otimes_{A}
\Gamma(\mfrak{U}, \mcal{M})$.
For such $\mfrak{U}$ one has
$\mrm{H}^{1}(\mfrak{U}, \mcal{M}) = 0$.
\end{cor}

\begin{proof}
Choose $\mfrak{U}$ affine such that $\mcal{M}|_{\mfrak{U}}$ has a
presentation
$\mcal{O}_{\mfrak{U}}^{(J)} \xrightarrow{\phi} \mcal{O}_{\mfrak{U}}^{(I)}
\xrightarrow{\psi} \mcal{M}|_{\mfrak{U}} \ar 0$. Define
$M := \opn{Coker}(\phi: A^{(I)} \ar A^{(J)})$.
Applying the exact functor $\mcal{O}_{\mfrak{U}} \otimes_{A}$ to
$A^{(I)} \xrightarrow{\phi} A^{(J)} \ar M \ar 0$
we get
$\mcal{M}|_{\mfrak{U}} \cong \mcal{O}_{\mfrak{U}} \otimes_{A} M$.
By the Proposition
$M \cong \Gamma(\mfrak{U}, \mcal{M})$.
As for $\mrm{H}^{1}(\mfrak{U}, - )$, use the fact that
it vanishes on coherent sheaves.
\end{proof}

\begin{prop} \label{prop3.3}
Let $\mcal{M}$ be coherent and $\mcal{N}$ quasi-coherent \textup{(}resp.\
coherent\textup{)}. Then
$\mcal{H}om_{\mcal{O}_{\mfrak{X}}}(\mcal{M}, \mcal{N})$
is quasi-coherent \textup{(}resp.\ coherent\textup{)}.
\end{prop}

\begin{proof}
For small enough $\mfrak{U} = \opn{Spf} A$ we get
$\mcal{M}|_{\mfrak{U}} \cong \mcal{O}_{\mfrak{U}} \otimes_{A} M$ and
$\mcal{N}|_{\mfrak{U}} \cong \mcal{O}_{\mfrak{U}} \otimes_{A} N$.
Now for any $\mfrak{V} = \opn{Spf} B \subset \mfrak{U}$, $A \ar B$ is flat;
so
\[ \opn{Hom}_{B}(B \otimes_{A} M, B \otimes_{A} N) \cong
B \otimes_{A} \opn{Hom}_{A}(M, N) . \]
Hence
\[ \mcal{H}om_{\mcal{O}_{\mfrak{X}}}(\mcal{M}, \mcal{N})|_{\mfrak{U}} \cong
\mcal{O}_{\mfrak{U}} \otimes_{A} \opn{Hom}_{A}(M,N) . \]
\end{proof}

Recall that a subcategory $\mathsf{B}$ of an abelian category $\mathsf{A}$
is called a thick abelian subcategory if for any exact sequence
$M_{1} \ar M_{2} \ar N \ar M_{3} \ar M_{4}$
in $\mathsf{A}$ with $M_{i} \in \mathsf{B}$, also $N \in \mathsf{B}$.

\begin{prop} \label{prop3.1}
The category $\mathsf{QCo}(\mfrak{X})$ is a thick abelian subcategory of
$\mathsf{Mod}(\mfrak{X})$.
\end{prop}

\begin{proof}
First observe that the kernel and cokernel of a homomorphism
$\mcal{M} \ar \mcal{N}$ between quasi-coherent sheaves is also
quasi-coherent.
This is immediate from Cor.\ \ref{cor3.1} and Prop.\ \ref{prop3.2}.
So it suffices to prove:
$0 \ar \mcal{M}' \ar \mcal{M} \ar \mcal{M}'' \ar 0$ exact,
$\mcal{M}',  \mcal{M}''$ quasi-coherent $\Rightarrow$ $\mcal{M}$
quasi-coherent. For a sufficiently small affine open formal subscheme
$\mfrak{U} = \opn{Spf} A$ we will get, by Cor.\ \ref{cor3.1}, that
$\mrm{H}^{1}(\mfrak{U}, \mcal{M}') = 0$. Hence the
sequence
\[ 0 \ar \Gamma(\mfrak{U}, \mcal{M}') \ar
M= \Gamma(\mfrak{U}, \mcal{M}) \ar
\Gamma(\mfrak{U}, \mcal{M}'') \ar 0 \]
is exact. This implies that
$\mcal{M}|_{\mfrak{U}} \cong \mcal{O}_{\mfrak{U}} \otimes_{A} M$.
\end{proof}

\begin{dfn} \label{dfn3.1}
Let $\mcal{M}$ be an $\mcal{O}_{\mfrak{X}}$-module. Define
\[ \ul{\Gamma}_{\mrm{disc}} \mcal{M} := \lim_{n \ar}
\mcal{H}om_{\mcal{O}_{\mfrak{X}}}(\mcal{O}_{\mfrak{X}} / \mcal{I}^{n},
\mcal{M}) \subset \mcal{M} \]
where $\mcal{I} \subset \mcal{O}_{\mfrak{X}}$ is any defining ideal.
$\mcal{M}$ is called {\em discrete} if
$\ul{\Gamma}_{\mrm{disc}} \mcal{M} = \mcal{M}$.
\end{dfn}

\begin{prop} \label{prop3.4}
Let $\mcal{M}$ be a quasi-coherent
$\mcal{O}_{\mfrak{X}}$-module. Then
$\ul{\Gamma}_{\mrm{disc}} \mcal{M}$ is quasi-coherent, and in fact
is a direct limit of discrete coherent $\mcal{O}_{\mfrak{X}}$-modules.
\end{prop}

\begin{proof}
Let $X_{n}$ be as in formula (\ref{eqn1.1}) and
$\mcal{M}_{n} := \mcal{H}om_{\mcal{O}_{\mfrak{X}}}(\mcal{O}_{X_{n}},
\mcal{M})$, so
$\ul{\Gamma}_{\mrm{disc}} \mcal{M} = \lim_{n \ar} \mcal{M}_{n}$.
If $\mcal{M}$ is quasi-coherent, then $\mcal{M}_{n}$ is a quasi-coherent
$\mcal{O}_{X_{n}}$-module (by Prop.\ \ref{prop3.3}), and hence is a
direct limit of coherent modules.
\end{proof}

% ** section 4 **

\section{Some Derived Functors of $\mcal{O}_{\mfrak{X}}$-Modules}

Denote by $\mathsf{Mod}_{\mrm{disc}}(\mfrak{X})$
(resp.\ $\mathsf{QCo}_{\mrm{disc}}(\mfrak{X})$)
the full subcategory of $\mathsf{Mod}(\mfrak{X})$ consisting
of discrete modules (resp.\ discrete quasi-coherent modules).
These are thick abelian subcategories.
In this section we study injective objects in the category
$\msf{QCo}_{\mrm{disc}}(\mfrak{X})$, and introduce the discrete Cousin
functor $\mrm{E} \mrm{R} \ul{\Gamma}_{\mrm{disc}}$.

\begin{lem} \label{lem4.1}
$\msf{Mod}_{\mrm{disc}}(\mfrak{X})$
is a locally noetherian category, with enough injectives.
\end{lem}

\begin{proof}
A family of noetherian generators consists of the sheaves
$\mcal{O}_{U}$, where $X \subset \mfrak{X}$ is a closed subscheme,
$U \subset X$ is an open set, and $\mcal{O}_{U}$ is extended by $0$
to all of $X$ (cf.\ \cite{RD} Theorem II.7.8).
If $\mcal{J} \in \msf{Mod}(\mfrak{X})$ is injective then
$\ul{\Gamma}_{\mrm{disc}} \mcal{J}$ is injective in
$\msf{Mod}_{\mrm{disc}}(\mfrak{X})$.
\end{proof}

Given a point $x \in \mfrak{X}$ let $J(x)$ be an injective hull
of the residue field $k(x)$ over the local ring
$\mcal{O}_{\mfrak{X},x}$,
and let $\mcal{J}(x)$ be the corresponding
$\mcal{O}_{\mfrak{X}}$-module. Then $\mcal{J}(x)$ is a
discrete quasi-coherent sheaf, constant on $\overline{\{ x \}}$,
and it is injective in $\msf{Mod}(\mfrak{X})$.

\begin{prop} \label{prop4.1}
\begin{enumerate}
\item
$\msf{QCo}_{\mrm{disc}}(\mfrak{X})$ is a locally noetherian
category with enough injectives.
\item Let $\mcal{J} \in \msf{QCo}_{\mrm{disc}}(\mfrak{X})$
be an injective object. Then $\mcal{J}$ is injective in
$\msf{Mod}_{\mrm{disc}}(\mfrak{X})$ and injective on
$\msf{Coh}(\mfrak{X})$.
For any
$\mcal{M} \in \msf{Mod}_{\mrm{disc}}(\mfrak{X})$
or
$\mcal{M} \in \msf{Coh}(\mfrak{X})$ the sheaf
$\mcal{H}om_{\mfrak{X}}(\mcal{M}, \mcal{J})$ is flasque.
\end{enumerate}
\end{prop}

\begin{proof}
1.\
Let $\mcal{N} \in \msf{QCo}_{\mrm{disc}}(\mfrak{X})$.
Choose a defining ideal $\mcal{I}$ of $\mfrak{X}$ and let
$X_{0}$ be the scheme $(\mfrak{X},
\mcal{O}_{\mfrak{X}} / \mcal{I})$.
Define $\mcal{N}_{0} :=
\mcal{H}om_{\mfrak{X}}(\mcal{O}_{X_{0}}, \mcal{N})$,
which is a quasi-coherent $\mcal{O}_{X_{0}}$-module.
Then the injective hull of $\mcal{N}_{0}$ in $\msf{Mod}(X_{0})$
is isomorphic to
$\bigoplus_{\alpha} \mcal{J}_{0}(x_{\alpha})$
for some $x_{\alpha} \in X_{0}$.
According to Proposition \ref{prop3.4},
$\msf{QCo}_{\mrm{disc}}(\mfrak{X})$
is locally noetherian, and this implies that
$\bigoplus_{\alpha} \mcal{J}(x_{\alpha})$
is an injective object in it. Now
$\mcal{N}_{0} \subset \mcal{N}$
and
$\mcal{N}_{0} \subset \bigoplus_{\alpha} \mcal{J}(x_{\alpha})$
are essential submodules, so there is some homomorphism
$\mcal{N} \ar \bigoplus_{\alpha} \mcal{J}(x_{\alpha})$,
which is necessarily injective and essential.

\noindent 2.\
If $\mcal{N} = \mcal{J}$ is injective in
$\msf{QCo}_{\mrm{disc}}(\mfrak{X})$, it follows that
$\mcal{J} \ar \bigoplus_{\alpha} \mcal{J}(x_{\alpha})$
is an isomorphism.
Since $\msf{Mod}_{\mrm{disc}}(\mfrak{X})$ is locally noetherian it
follows that $\mcal{J}$ is injective in it.
Given $\mcal{M} \in \msf{Mod}_{\mrm{disc}}(\mfrak{X})$
and open sets
$\mfrak{V} \subset \mfrak{U} \subset \mfrak{X}$
consider the sheaves
$\mcal{M}|_{\mfrak{V}} \subset \mcal{M}|_{\mfrak{U}} \subset
\mcal{M}$
(extension by $0$). Then
$\mcal{H}om_{\mfrak{X}}(\mcal{M}|_{\mfrak{U}}, \mcal{J}) \ar
\mcal{H}om_{\mfrak{X}}(\mcal{M}|_{\mfrak{V}}, \mcal{J})$
is surjective.

The category
$\msf{Coh}(\mfrak{X})$ is noetherian, and therefore the functor
$\opn{Hom}_{\mfrak{X}}(-, \mcal{J})$ is exact on it.
Given $\mcal{M} \in \msf{Coh}(\mfrak{X})$ we have
$\mcal{H}om_{\mfrak{X}}(\mcal{M}, \mcal{J}) \cong
\bigoplus \mcal{H}om_{\mfrak{X}}(\mcal{M}, \mcal{J}(x_{\alpha}))$
which is clearly flasque.
\end{proof}

\begin{cor} \label{cor4.1}
Let
$\mcal{J}^{\bdot} \in \msf{D}^{+}(\msf{QCo}_{\mrm{disc}}(\mfrak{X}))$
be a complex of injectives. Then for any
$\mcal{M}^{\bdot} \in \msf{Mod}_{\mrm{disc}}(\mfrak{X})$
or
$\mcal{M}^{\bdot} \in \msf{Coh}(\mfrak{X})$ one has
\[ \begin{aligned}
\mrm{R} \mcal{H}om_{\mfrak{X}}(\mcal{M}^{\bdot}, \mcal{J}^{\bdot})
& \cong \mcal{H}om_{\mfrak{X}}(\mcal{M}^{\bdot}, \mcal{J}^{\bdot}) \\
\mrm{R} \opn{Hom}_{\mfrak{X}}(\mcal{M}^{\bdot}, \mcal{J}^{\bdot})
& \cong \opn{Hom}_{\mfrak{X}}(\mcal{M}^{\bdot}, \mcal{J}^{\bdot})
\cong \Gamma(\mfrak{X},
\mcal{H}om_{\mfrak{X}}(\mcal{M}^{\bdot}, \mcal{J}^{\bdot})) .
\end{aligned} \]
\end{cor}

\begin{proof}
The first equality follows from Proposition \ref{prop4.1}
(cf.\ \cite{RD} Section I.6). Since each sheaf
$\mcal{H}om_{\mfrak{X}}(\mcal{M}^{p}, \mcal{J}^{q})$
is flasque we obtain the second equality.
\end{proof}

The functor
$\ul{\Gamma}_{\mrm{disc}} : \msf{Mod}(\mfrak{X}) \ar
\msf{Mod}_{\mrm{disc}}(\mfrak{X})$
has a derived functor
\[ \mrm{R} \ul{\Gamma}_{\mrm{disc}} : \msf{D}^{+}(\msf{Mod}(\mfrak{X}))
\ar \msf{D}^{+}(\msf{Mod}_{\mrm{disc}}(\mfrak{X})) , \]
which is calculated by injective resolutions.

There is another way to compute cohomology with supports.
Let $t$ be an indeterminate. Define
$\mbf{K}^{\bdot}(t)$ to be the Koszul complex
$\mbb{Z}\sqbr{t} \xrightarrow{\cdot t} \mbb{Z}\sqbr{t}$, in dimensions $0$ and $1$,
and let
$\mbf{K}^{\bdot}_{\infty}(t) := \lim_{i \ar} \mbf{K}^{\bdot}(t^{i})$.
Given a sequence $\ul{t} = (t_{1}, \ldots, t_{n})$
define
$\mbf{K}^{\bdot}_{\infty}(\ul{t}) := \mbf{K}^{\bdot}_{\infty}(t_{1})
\otimes \cdots \otimes \mbf{K}^{\bdot}_{\infty}(t_{n})$,
a complex of flat $\mbb{Z} \sqbr{\ul{t}}$-modules
(in fact it's a commutative DGA).
If $A$ is a noetherian commutative ring and
$\ul{a} = (a_{1}, \ldots, a_{n}) \in A^{n}$, then we write
$\mbf{K}^{\bdot}_{\infty}(\ul{a})$ instead of
$\mbf{K}^{\bdot}_{\infty}(\ul{t}) \otimes_{\mbb{Z} \sqbr{\ul{t}}} A$.
Now suppose $\mfrak{a} \subset A$ is an ideal, and $\ul{a}$ are
generators of $\mfrak{a}$. Then for any
$M^{\bdot} \in \msf{D}^{+}(\msf{Mod}(A))$
there is a natural isomorphism
\begin{equation} \label{eqn4.1}
\mrm{R} \Gamma_{\mfrak{a}} M^{\bdot} \cong
\mbf{K}^{\bdot}_{\infty}(\ul{a}) \otimes M^{\bdot}
\end{equation}
in $\msf{D}(\msf{Mod}(A))$.
We refer to \cite{LS1}, \cite{Hg1} and \cite{AJL1} for full details and
proofs. For sheaves one has:

\begin{lem} \label{lem4.3}
Suppose $\ul{a} \in \Gamma(\mfrak{U}, \mcal{O}_{\mfrak{U}})^{n}$
generates a defining ideal of the formal sche\-me $\mfrak{U}$.
Then for any $\mcal{M}^{\bdot} \in \msf{D}^{+}(\msf{Mod}(\mfrak{U}))$
there is a natural isomorphism
\[ \mrm{R} \ul{\Gamma}_{\mrm{disc}} \mcal{M}^{\bdot} \cong
\mbf{K}^{\bdot}_{\infty}(\ul{a}) \otimes \mcal{M}^{\bdot} . \]
\end{lem}

\begin{proof}
Let $\mcal{I} := \mcal{O}_{\mfrak{U}} \cdot \ul{a}$. Then
$\ul{\Gamma}_{\mrm{disc}} = \ul{\Gamma}_{\mcal{I}}$, and we may use
\cite{AJL1} Lemma 3.1.1.
\end{proof}

\begin{prop} \label{prop4.2}
Let $X$ be a noetherian scheme, $X_{0} \subset X$ a closed subset,
$\mfrak{X} = X_{/ X_{0}}$ and
$g : \mfrak{X} \ar X$ the completion morphism. Then for any
$\mcal{M}^{\bdot} \in \msf{D}^{+}_{\mrm{qc}}(\msf{Mod}(X))$
there is a natural isomorphism
$g^{*} \mrm{R} \ul{\Gamma}_{X_{0}} \mcal{M}^{\bdot} \cong
\mrm{R} \ul{\Gamma}_{\mrm{disc}} g^{*} \mcal{M}^{\bdot}$.
In particular for a single quasi-coherent sheaf $\mcal{M}$ one has
$g^{*} \ul{\Gamma}_{X_{0}} \mcal{M} \cong
\ul{\Gamma}_{\mrm{disc}} g^{*} \mcal{M}$.
\end{prop}

\begin{proof}
Let
$\mcal{M}^{\bdot} \ar \mcal{J}^{\bdot}$ be a resolution by quasi-coherent
injectives. Since $g$ is flat we get
\[ \phi :
g^{*} \mrm{R} \ul{\Gamma}_{X_{0}} \mcal{M}^{\bdot} =
g^{*} \ul{\Gamma}_{X_{0}} \mcal{J}^{\bdot} \ar
\ul{\Gamma}_{\mrm{disc}} g^{*} \mcal{J}^{\bdot} \ar
\mrm{R} \ul{\Gamma}_{\mrm{disc}} g^{*} \mcal{J}^{\bdot}
= \mrm{R} \ul{\Gamma}_{\mrm{disc}} g^{*} \mcal{M}^{\bdot} . \]
Locally on any affine open $U \subset X$, with $U_{0} = U \cap X_{0}$
and $\mfrak{U} = U_{/U_{0}}$, we can find $\ul{a}$ in
$\Gamma(U, \mcal{O}_{U})$ which define $U_{0}$. It's known that
$\ul{\Gamma}_{U_{0}} (\mcal{J}^{\bdot}|_{U}) \ar
\mbf{K}^{\bdot}_{\infty}(\ul{a}) \otimes (\mcal{J}^{\bdot}|_{U})$
is a quasi-isomorphism. Since $g$ is flat we obtain quasi-isomorphisms
\begin{multline*}
\hspace{10mm}
\phi|_{\mfrak{U}} : g^{*} \ul{\Gamma}_{U_{0}} (\mcal{J}^{\bdot}|_{U})
\ar g^{*} \left( \mbf{K}^{\bdot}_{\infty}(\ul{a}) \otimes
(\mcal{J}^{\bdot}|_{U}) \right) \\
\cong \mbf{K}^{\bdot}_{\infty}(\ul{a}) \otimes
g^{*} (\mcal{J}^{\bdot}|_{U}) =
\mrm{R} \ul{\Gamma}_{\mrm{disc}} g^{*} (\mcal{J}^{\bdot}|_{U}) .
\hspace{10mm}
\end{multline*}
It follows that $\phi$ is an isomorphism.
\end{proof}

Denote by $\msf{D}^{+}_{\mrm{d}}(\msf{Mod}(\mfrak{X}))$ the
subcategory of complexes with discrete cohomologies.

\begin{lem} \label{lem4.2}
\begin{enumerate}
\item If
$\mcal{M}^{\bdot} \in \msf{D}^{+}_{\mrm{d}}(\msf{Mod}(\mfrak{X}))$
then
$\mrm{R} \ul{\Gamma}_{\mrm{disc}} \mcal{M} \ar \mcal{M}$
is an isomorphism.
\item If
$\mcal{M}^{\bdot} \in \msf{D}^{+}_{\mrm{qc}}(\msf{Mod}(\mfrak{X}))$
then
$\mrm{R} \ul{\Gamma}_{\mrm{disc}} \mcal{M} \in
\msf{D}^{+}_{\mrm{qc}}(\msf{Mod}_{\mrm{disc}}(\mfrak{X}))$.
\end{enumerate}
\end{lem}

\begin{proof}
From Lemma \ref{lem4.3} we see that the functor
$\mrm{R} \ul{\Gamma}_{\mrm{disc}}$ has finite cohomological dimension.
By way-out reasons (cf.\ \cite{RD} Section I.7) we may assume
$\mcal{M}^{\bdot}$ is a single discrete (resp.\ quasi-coherent) sheaf.
Then the claims are obvious (use Proposition \ref{prop3.4} for 2).
\end{proof}

\begin{thm} \label{thm4.1}
The identity functor
$\msf{D}^{+}(\msf{QCo}_{\mrm{disc}}(\mfrak{X})) \ar
\msf{D}^{+}_{\mrm{dqc}}(\msf{Mod}(\mfrak{X}))$
is an equivalence of categories.
In particular any
$\mcal{M}^{\bdot} \in \msf{D}^{+}_{\mrm{dqc}}(\msf{Mod}(\mfrak{X}))$
is isomorphic to a complex of injectives
$\mcal{J}^{\bdot} \in \msf{D}^{+}(\msf{QCo}_{\mrm{disc}}(\mfrak{X}))$.
\end{thm}

\begin{proof}
According to Lemma \ref{lem4.2} we see that
$\msf{D}^{+}_{\mrm{qc}}(\msf{Mod}_{\mrm{disc}}(\mfrak{X})) \ar
\msf{D}^{+}_{\mrm{dqc}}(\msf{Mod}(\mfrak{X}))$
is an equivalence with quasi-inverse $\mrm{R} \ul{\Gamma}_{\mrm{disc}}$.
Next, by Proposition \ref{prop4.1} and by \cite{RD} Proposition I.4.8,
the functor
$\msf{D}^{+}(\msf{QCo}_{\mrm{disc}}(\mfrak{X})) \ar
\msf{D}^{+}_{\mrm{qc}}(\msf{Mod}_{\mrm{disc}}(\mfrak{X}))$
is an equivalence.
\end{proof}

\begin{rem} \label{rem6.10}
In \cite{AJL2} it is proved that
$\msf{D}(\msf{QCo}_{\mrm{disc}}(\mfrak{X}))
\ar \msf{D}_{\mrm{dqc}}(\msf{Mod}(\mfrak{X}))$
is an equivalence, using the quasi-coherator functor.
\end{rem}

Suppose there is a codimension function
$d : \mfrak{X} \ar \mbb{Z}$,
i.e.\ a function satisfying $d(y) = d(x) + 1$ whenever $(x, y)$ is
an immediate specialization pair. Then there is a filtration
$\cdots \supset Z^{p} \supset Z^{p+1} \supset \cdots$
of $\mfrak{X}$, with
$Z^{p} := \{ F \subset \mfrak{X} \mid F \text{ closed}, d(F) \geq p \}$.
Here $d(F) := \opn{min} \{ d(x) \mid x \in F \}$.
This filtration determines a Cousin functor
\begin{equation}
\mrm{E} : \msf{D}^{+}(\msf{Ab}(\mfrak{X})) \ar
\msf{C}^{+}(\msf{Ab}(\mfrak{X}))
\end{equation}
where $\msf{C}^{+}$ denotes the abelian category of bounded below
complexes (cf.\ \cite{RD} \S IV.1).

Given a point $x \in \mfrak{X}$ and a sheaf
$\mcal{M} \in \msf{Ab}(\mfrak{X})$ we let
$\Gamma_{x} \mcal{M} := (\ul{\Gamma}_{\, \overline{\{x\}}\, }
\mcal{M})_{x}$ $\subset \mcal{M}_{x}$.
The derived functor
$\mrm{R} \Gamma_{x} : \msf{D}^{+}(\msf{Ab}(\mfrak{X})) \ar
\msf{D}(\msf{Ab})$
is calculated by flasque sheaves. Let us write
$\mrm{H}_{x}^{q} \mcal{M} :=
\mrm{H}^{q} \mrm{R} \Gamma_{x} \mcal{M}$,
the local cohomology, and let $i_{x} : \{x\} \ar \mfrak{X}$
be the inclusion

According to \cite{RD} \S IV.1 Motif F one has a natural isomorphism
\begin{equation} \label{eqn4.3}
\mrm{E}^{p} \mcal{M}^{\bdot} =
\mcal{H}_{Z^{p} / Z^{p+1}}^{p} \mcal{M}^{\bdot} \cong
\bigoplus_{d(x) = p} i_{x *} \mrm{H}_{x}^{p} \mcal{M}^{\bdot} .
\end{equation}

Observe that if
$\mcal{M} \in \msf{D}^{+}(\msf{Mod}(\mfrak{X}))$
then
$\mrm{E} \mcal{M}^{\bdot} \in
\msf{C}^{+}(\msf{Mod}(\mfrak{X}))$
and
$\mrm{R} \Gamma_{x} \mcal{M} \in$ \newline
$\msf{D}^{+}(\msf{Mod}(\mcal{O}_{\mfrak{X}, x}))$.

Unlike an ordinary scheme, on a formal scheme the topological
support of a quasi-coherent sheaf does not coincide with its
algebraic support. But for discrete sheaves these two notions of support
do coincide. This suggests:

\begin{dfn}
Given $\mcal{M} \in \msf{D}^{+}(\msf{Mod}(\mfrak{X}))$
its {\em discrete Cousin complex} is \newline
$\mrm{E} \mrm{R} \ul{\Gamma}_{\mrm{disc}} \mcal{M}^{\bdot}$.
\end{dfn}

\begin{thm} \label{thm4.2}
For any
$\mcal{M}^{\bdot} \in \msf{D}^{+}_{\mrm{qc}}(\msf{Mod}(\mfrak{X}))$
the complex
$\mrm{E} \mrm{R} \ul{\Gamma}_{\mrm{disc}} \mcal{M}^{\bdot}$
consists of discrete quasi-coherent sheaves. So we get a functor
\[ \mrm{E} \mrm{R} \ul{\Gamma}_{\mrm{disc}} :
\msf{D}^{+}_{\mrm{qc}}(\msf{Mod}(\mfrak{X})) \ar
\msf{C}^{+}(\msf{QCo}_{\mrm{disc}}(\mfrak{X})). \]
\end{thm}

\begin{proof}
According to Theorem \ref{thm4.1} we may assume
$\mcal{N}^{\bdot} = \mrm{R} \ul{\Gamma}_{\mrm{disc}} \mcal{M}^{\bdot}$
is in \newline
$\msf{D}^{+}(\msf{QCo}_{\mrm{disc}}(\mfrak{X}))$.
On any open formal subscheme $\mfrak{U} = \opn{Spf} A$ we get
$\mcal{N}^{\bdot} = \mcal{O}_{\mfrak{U}} \otimes_{A} N^{\bdot}$,
where
$N^{q} = \Gamma(\mfrak{U}, N^{q})$
(cf.\ Propositions \ref{prop3.4} and \ref{prop3.2})
Then for $x \in \mfrak{U}$,
\[ \mrm{R} \Gamma_{x} \mrm{R} \ul{\Gamma}_{\mrm{disc}} \mcal{M}^{\bdot}
= \mrm{R} \Gamma_{x} \mcal{N}^{\bdot} =
\mrm{R} \Gamma_{\mfrak{p}} N^{\bdot}_{\mfrak{p}} \]
where $\mfrak{p} \subset A$ is the prime ideal of $x$. Hence
$\mrm{H}^{q}_{x} \mrm{R} \ul{\Gamma}_{\mrm{disc}} \mcal{M}^{\bdot}
= \mrm{H}^{q}_{\mfrak{p}} N^{\bdot}_{\mfrak{p}}$
is $\mfrak{p}$-torsion. So the sheaf corresponding to $x$ in
(\ref{eqn4.3}) is quasi-coherent and discrete.
\end{proof}

% ** section 5 **

\section{Dualizing Complexes on Formal Schemes}

In this section we propose a theory of duality on noetherian formal
sche\-mes. There is a fundamental difference between this theory and the
duality theory on schemes, as developed in \cite{RD}. A dualizing
complex $\mcal{R}^{\bdot}$ on a scheme $X$ has coherent cohomology
sheaves; this will not be true on a general formal scheme $\mfrak{X}$,
where $\mrm{H}^{q} \mcal{R}^{\bdot}$ are discrete quasi-coherent sheaves
(Def.\ \ref{dfn5.1}).
We prove uniqueness of dualizing complexes (Thm.\ \ref{thm5.1}),
and existence in some cases (Prop.\ \ref{prop5.8} and Thm.\ \ref{thm5.3}).

Before we begin here is an instructive example due to J.\ Lipman.

\begin{exa} \label{exa5.2}
Consider the ring $A = k[\sqbr{t}]$ of formal power series over a field
$k$. Let $\mfrak{X} := \opn{Spf} A$, which has a single point. The modules
$A$ and $J = \mrm{H}^{1}_{(t)} A$ both have finite injective dimension
and satisfy $\opn{Hom}_{A}(A, A) = \opn{Hom}_{A}(J, J) = A$. Which
one is a dualizing complex on $\mfrak{X}$? We will see that $J$ is
the correct answer (Def.\ \ref{dfn5.1}), and $A$ is a ``fake'' dualizing
complex (Thm.\ \ref{thm5.3}). The relevant relation between them is:
$J = \mrm{R} \Gamma_{\mrm{disc}} A [1]$.
\end{exa}

Suppose
$\mcal{N}^{\bdot} \in
\msf{D}^{+}(\msf{Mod}_{\mrm{disc}}(\mfrak{X}))$.
We say $\mcal{N}^{\bdot}$ has finite injective dimension on
$\msf{QCo}_{\mrm{disc}}(\mfrak{X})$
if there is an integer $q_{0}$ s.t.\ for all $q > q_{0}$ and
$\mcal{M} \in  \msf{QCo}_{\mrm{disc}}(\mfrak{X})$,
$\mrm{H}^{q} \mrm{R} \opn{Hom}_{\mfrak{X}}
(\mcal{M}, \mcal{N}^{\bdot}) = 0$.

\begin{dfn} \label{dfn5.1}
A {\em dualizing complex} on $\mfrak{X}$ is a complex
$\mcal{R}^{\bdot} \in
\msf{D}^{\mrm{b}}_{\mrm{dqc}}(\msf{Mod}(\mfrak{X}))$
satisfying:
\begin{enumerate}
\rmitem{i} $\mcal{R}^{\bdot}$ has finite injective dimension
on $\msf{QCo}_{\mrm{disc}}(\mfrak{X})$.
\rmitem{ii} The adjunction morphism
$\mcal{O}_{\mfrak{X}} \ar
\mrm{R} \mcal{H}om_{\mfrak{X}}
(\mcal{R}^{\bdot}, \mcal{R}^{\bdot})$
is an isomorphism.
\rmitem{iii} For some defining ideal $\mcal{I}$ of $\mfrak{X}$,
$\mrm{R} \mcal{H}om_{\mfrak{X}}(\mcal{O}_{\mfrak{X}} / \mcal{I},
\mcal{R}^{\bdot})$
has coherent cohomology sheaves.
\end{enumerate}
\end{dfn}

\begin{lem} \label{lem5.2}
Let
$\mcal{N}^{\bdot} \in
\msf{D}^{+}_{\mrm{dqc}}(\msf{Mod}(\mfrak{X}))$.
Then $\mcal{N}^{\bdot}$ has finite injective dimension on
$\msf{QCo}_{\mrm{disc}}(\mfrak{X})$ iff it is isomorphic to a bounded
complex of injectives in $\msf{QCo}_{\mrm{disc}}(\mfrak{X})$.
\end{lem}

\begin{proof}
Because of Theorem \ref{thm4.1}, the proof is just like
\cite{RD} Prop.\ I.7.6.
\end{proof}

In light of this, we can, when convenient, assume the dualizing complex
$\mcal{R}^{\bdot}$ is a bounded complex of discrete quasi-coherent
injectives.

\begin{prop} \label{prop5.2}
Let $\mcal{R}^{\bdot}$ be a dualizing complex on $\mfrak{X}$.
Then for any
$\mcal{M}^{\bdot} \in \msf{D}^{\mrm{b}}_{\mrm{c}}(\msf{Mod}(\mfrak{X}))$
the morphism of adjunction
\[ \mcal{M}^{\bdot} \ar
\mrm{R} \mcal{H}om_{\mfrak{X}} (
\mrm{R} \mcal{H}om_{\mfrak{X}} (\mcal{M}^{\bdot},
\mcal{R}^{\bdot}), \mcal{R}^{\bdot}) \]
is an isomorphism.
\end{prop}

\begin{proof}
We can assume $\mfrak{X}$ is affine, and so replace $\mcal{M}^{\bdot}$
with a complex of coherent sheaves.
By ``way-out'' arguments (cf.\ \cite{RD} Section I.7)
we reduce to the case $\mcal{M}^{\bdot} = \mcal{O}_{\mfrak{X}}$,
which property (ii) applies.
\end{proof}

\begin{lem} \label{lem5.4}
Suppose $\mcal{R}^{\bdot}$ is a dualizing complex on $\mfrak{X}$.
Let $\mcal{I}$ be any defining ideal of $\mfrak{X},$ and let $X_{0}$ be
the scheme
$(\mfrak{X}, \mcal{O}_{\mfrak{X}} / \mcal{I})$.
Then
$\mrm{R} \mcal{H}om_{\mfrak{X}}
(\mcal{O}_{X_{0}}, \mcal{R}^{\bdot})$
is a dualizing complex on $X_{0}$.
\end{lem}

\begin{proof}
We can assume $\mcal{R}^{\bdot}$ is a bounded complex of injectives
in $\msf{QCo}_{\mrm{disc}}(\mfrak{X})$, so
$\mcal{R}^{\bdot}_{0} := \mcal{H}om_{\mfrak{X}}
(\mcal{O}_{X_{0}}, \mcal{R}^{\bdot})$
is a complex of injectives on $X_{0}$. Property (iii) implies that
$\mcal{R}^{\bdot}_{0}$ has coherent cohomology sheaves. Now
\[ \mcal{H}om_{X_{0}}(\mcal{R}^{\bdot}_{0}, \mcal{R}^{\bdot}_{0})
\cong \mcal{H}om_{\mfrak{X}} (
\mcal{H}om_{\mfrak{X}} (\mcal{O}_{X_{0}},
\mcal{R}^{\bdot}), \mcal{R}^{\bdot})
\cong \mcal{O}_{X_{0}} , \]
so $\mcal{R}^{\bdot}_{0}$ is dualizing.
\end{proof}

\begin{thm} \label{thm5.1} \textup{(Uniqueness)}\
Suppose $\mcal{R}^{\bdot}$ and
$\tilde{\mcal{R}}^{\bdot}$ are dualizing complexes and $\mfrak{X}$
is connected. Then
$\tilde{\mcal{R}}^{\bdot} \cong \mcal{R}^{\bdot} \otimes
\mcal{L}[n]$
in
$\msf{D}(\msf{Mod}(\mfrak{X}))$,
for some invertible sheaf $\mcal{L}$ and integer $n$.
\end{thm}

\begin{proof}
We can assume both $\mcal{R}^{\bdot}$ and $\tilde{\mcal{R}}^{\bdot}$
are bounded complexes of injectives in
$\msf{QCo}_{\mrm{disc}}(\mfrak{X})$.
Choose a defining ideal $\mcal{I}$ and let
$X_{m}$ be the scheme
$(\mfrak{X},  \mcal{O}_{\mfrak{X}} / \mcal{I}^{m+1})$.
Define a complex
$\mcal{R}^{\bdot}_{m} := \mcal{H}om_{\mfrak{X}}
(\mcal{O}_{X_{m}}, \mcal{R}^{\bdot})$
and likewise $\tilde{\mcal{R}}^{\bdot}_{m}$.
These are dualizing complexes on $X_{m}$, so by \cite{RD} Thm.\ IV.3.1
there is an isomorphism
\[ \phi_{m} : \mcal{R}^{\bdot}_{m} \otimes \mcal{L}_{m}[n_{m}]
\ar  \tilde{\mcal{R}}^{\bdot}_{m} \]
in $\msf{D}(\msf{Mod}(X_{m}))$, for some
invertible sheaf $\mcal{L}_{m}$ and integer $n_{m}$.
Writing
$\mcal{M}_{m}^{\bdot} := \mcal{H}om_{X_{m}}(
\mcal{R}^{\bdot}_{m}, \tilde{\mcal{R}}^{\bdot}_{m})$
we have
$\mcal{M}_{m}^{\bdot} \cong \mcal{L}_{m}[n_{m}]$
in $\msf{D}(\msf{Mod}(X_{m}))$.
Now
\[ \mcal{M}_{m}^{\bdot} \cong
\mcal{H}om_{X_{m+1}}(
\mcal{H}om_{X_{m+1}}(\mcal{O}_{X_{m}},
\mcal{R}^{\bdot}_{m+1}), \mcal{R}^{\bdot}_{m+1})) \otimes
\mcal{L}_{m+1}[n_{m+1}] \]
as complexes of $\mcal{O}_{X_{m+1}}$-modules, so by the dualizing
property of $\mcal{R}^{\bdot}_{m+1}$  we deduce an isomorphism
$\mcal{M}_{m}^{\bdot} \cong \mcal{O}_{X_{m}}
\otimes \mcal{L}_{m+1}[n_{m+1}]$
in $\msf{D}(\msf{Mod}(X_{m+1}))$.
We conclude that
$n_{m} = n_{m+1}$ and
$\mcal{L}_{m} \cong \mcal{O}_{X_{m}} \otimes \mcal{L}_{m+1}$.
Set $n := n_{m}$  and
$\mcal{L} := \lim_{\leftarrow m} \mcal{L}_{m}$.

Next, since
$\mcal{R}^{q}_{m} \subset \mcal{R}^{q}_{m+1}$
and $\tilde{\mcal{R}}^{q}_{m+1}$
is injective in $\msf{Mod}(X_{m+1})$, we see that
$\mcal{M}_{m+1}^{q} \ar \mcal{M}_{m}^{q}$
is surjective for all $q,m$. Furthermore,
$\mrm{H}^{q} \mcal{M}^{\bdot}_{m+1} \ar
\mrm{H}^{q} \mcal{M}^{\bdot}_{m}$
is also surjective, since
$\mrm{H}^{q} \mcal{M}^{\bdot}_{m} = \mcal{L}_{m}$ or $0$.
Define
\[ \mcal{M}^{\bdot} := \mcal{H}om_{\mfrak{X}}(
\mcal{R}^{\bdot}, \tilde{\mcal{R}}^{\bdot})
\cong \lim_{\leftarrow m} \mcal{M}^{\bdot}_{m} . \]
According to \cite{Ha} Cor.\ I.4.3 and Prop.\ I.4.4 it follows that
$\mrm{H}^{q} \mcal{M}^{\bdot} = \lim_{\leftarrow m}
\mrm{H}^{q} \mcal{M}^{\bdot}_{m}$.
This implies that
$\mcal{H}om_{\mfrak{X}}(
\mcal{R}^{\bdot} \otimes \mcal{L}[n], \tilde{\mcal{R}}^{\bdot}))
\cong \mcal{O}_{\mfrak{X}}$
in $\msf{D}(\msf{Mod}(\mfrak{X}))$, so by Corollary \ref{cor4.1}
\[ \mrm{H}^{0} \opn{Hom}_{\mfrak{X}}
(\mcal{R}^{\bdot}\otimes \mcal{L}[n], \tilde{\mcal{R}}^{\bdot})
\cong \Gamma(\mfrak{X}, \mcal{O}_{\mfrak{X}}) .
\]
Choose a homomorphism of complexes
$\phi : \mcal{R}^{\bdot} \otimes \mcal{L}[n] \ar
\tilde{\mcal{R}}^{\bdot}$
corresponding to
$1 \in \Gamma(\mfrak{X}, \mcal{O}_{\mfrak{X}})$.
Backtracking we see that for every $m$, $\phi$ induces a homomorphism
$\mcal{R}^{\bdot}_{m} \otimes \mcal{L}[n] \ar
\tilde{\mcal{R}}^{\bdot}_{m}$
which represents $\phi_{m}$ in
$\msf{D}(\msf{Mod}(X_{m}))$. So
$\phi = \lim_{m \ar} \phi_{m}$ is a quasi-isomorphism.
\end{proof}

\begin{prob}
Let $\mcal{R}^{\bdot}$ be a dualizing complex. Is it true that the following
conditions on
$\mcal{N}^{\bdot} \in \msf{D}^{\mrm{b}}_{\mrm{dqc}}(\msf{Mod}(\mfrak{X}))$
are equivalent?
\begin{enumerate}
\rmitem{i}
$\mcal{N}^{\bdot} \cong \mrm{R} \mcal{H}om_{\mfrak{X}}^{\bdot}(
\mcal{M}^{\bdot}, \mcal{R}^{\bdot})$
for some
$\mcal{M}^{\bdot} \in \msf{D}^{\mrm{b}}_{\mrm{c}}(\msf{Mod}(\mfrak{X}))$.
\rmitem{ii} For any $\mcal{M}$ discrete coherent,
$\mrm{R} \mcal{H}om_{\mfrak{X}}^{\bdot}(
\mcal{M}, \mcal{N}^{\bdot}) \in
\msf{D}^{\mrm{b}}_{\mrm{c}}(\msf{Mod}(\mfrak{X}))$.
\end{enumerate}
\end{prob}

Recall that for a point $x \in \mfrak{X}$ we denote by $J(x)$ an
injective hull of $k(x)$ over $\mcal{O}_{\mfrak{X}, x}$,
and $\mcal{J}(x)$ is the corresponding quasi-coherent sheaf.

\begin{lem} \label{lem5.6}
Suppose $\mcal{R}^{\bdot}$ is a dualizing complex on $\mfrak{X}$.
For any $x \in \mfrak{X}$ there is a unique integer $d(x)$ s.t.\
\[ \mrm{H}^{q}_{x} \mcal{R}^{\bdot} \cong
\begin{cases}
J(x) & \text{ if } q = d(x)\\
0 & \text{ otherwise}.
\end{cases} \]
Furthermore $d$ is a codimension function.
\end{lem}

\begin{proof}
We can assume $\mcal{R}^{\bdot}$ is a bounded complex of injectives
in $\msf{QCo}_{\mrm{disc}}(\mfrak{X})$. Then as seen before
$\mrm{H}^{q}_{x} \mcal{R}^{\bdot} =
\mrm{H}^{q} \Gamma_{x} \mcal{R}^{\bdot}$. Define schemes $X_{m}$
and complexes $\mcal{R}^{\bdot}_{m}$ like in the proof of Thm.\
\ref{thm5.1}. Since $\mcal{R}^{\bdot}_{m}$ is dualizing it determines
a codimension function $d_{m}$ on $X_{m}$ (cf.\ \cite{RD} Ch.\ V \S 7).
But the arguments used before show that $d_{m} = d_{m+1}$. Finally
$\mrm{H}^{q} \Gamma_{x} \mcal{R}^{\bdot} =
\lim_{m \ar} \mrm{H}^{q} \Gamma_{x} \mcal{R}^{\bdot}_{m}$,
and
$\mrm{H}^{q} \Gamma_{x} \mcal{R}^{\bdot}_{m} \cong J_{m}(x)$,
an injective hull of $k(x)$ over $\mcal{O}_{X_{m}, x}$.
\end{proof}

\begin{dfn} \label{dfn5.3}
A residual complex on the noetherian formal scheme $\mfrak{X}$
is a dualizing complex $\mcal{K}^{\bdot}$ which
is isomorphic, as $\mcal{O}_{\mfrak{X}}$-module, to
$\bigoplus_{x \in \mfrak{X}} \mcal{J}(x)$.
\end{dfn}

\begin{prop} \label{prop5.7}
Say $\mcal{R}^{\bdot}$ is a dualizing complex on $\mfrak{X}$.
Let $d$ be the codimension function above, and let $\mrm{E}$ be the
associated Cousin functor. Then
$\mcal{R}^{\bdot} \cong \mrm{E} \mcal{R}^{\bdot}$ in
$\msf{D}(\msf{Mod}(\mfrak{X}))$, and
$\mrm{E} \mcal{R}^{\bdot}$ is a residual complex.
\end{prop}

\begin{proof}
By Lemma \ref{lem5.6} $\mcal{R}^{\bdot}$ is a Cohen-Macaulay complex,
in the sense of \cite{RD} p.\ 247, Definition. So there exists some
isomorphism
$\mcal{R}^{\bdot} \ar \mrm{E} \mcal{R}^{\bdot}$
in $\msf{D}^{\mrm{b}}(\msf{Mod}(\mfrak{X}))$.
\end{proof}

To conclude this section we consider some situations where a dualizing
complex exists. If $f : \mfrak{X} \ar \mfrak{Y}$ is a morphism
then $(\mfrak{Y}, f_{*} \mcal{O}_{\mfrak{X}})$ is a ringed space,
and
$\bar{f} : \mfrak{X} \ar (\mfrak{Y}, f_{*} \mcal{O}_{\mfrak{X}})$
is a morphism of ringed spaces.

\begin{prop} \label{prop5.8}
Let $f : \mfrak{X} \ar \mfrak{Y}$ be a formally finite morphism,
and assume $\mcal{K}^{\bdot}$ is a residual complex on $\mfrak{Y}$.
Then
$\bar{f}^{*} \mcal{H}om_{\mfrak{Y}}(f_{*} \mcal{O}_{\mfrak{X}},
\mcal{K}^{\bdot})$
is a residual complex on $\mfrak{X}$.
\end{prop}

\begin{proof}
Let $f_{n} : X_{n} \ar Y_{n}$ be morphisms as in Lemma \ref{lem1.12},
and let
$\mcal{K}_{n}^{\bdot}$ \linebreak
$:= \mcal{H}om_{\mfrak{Y}}(\mcal{O}_{Y_{n}},
\mcal{K}^{\bdot})$.
Since $f_{n}$ is a finite morphism,
$\bar{f}_{n}^{*} \mcal{H}om_{Y_{n}}(f_{n *} \mcal{O}_{X_{n}},
\mcal{K}_{n}^{\bdot})$
is a residual complex on $X_{n}$. As in the proof of Thm.\ \ref{thm5.1},
\[ \bar{f}^{*} \mcal{H}om_{\mfrak{Y}}(f_{*} \mcal{O}_{\mfrak{X}},
\mcal{K}^{\bdot}) \cong \lim_{n \ar}
\bar{f}_{n}^{*} \mcal{H}om_{Y_{n}}(f_{n *} \mcal{O}_{X_{n}},
\mcal{K}_{n}^{\bdot}) \]
is residual.
\end{proof}

\begin{exa} \label{exa5.1}
Suppose $X_{0} \subset X$ is closed, $\mfrak{X} = X_{/ X_{0}}$
and $g : \mfrak{X} \ar X$ is the completion morphism.
Let $\mcal{K}^{\bdot}$ be a residual complex on $X$. In this case
$g = \bar{g}$, and by Proposition \ref{prop4.2}
\[ g^{*} \mcal{H}om_{X}(g_{*} \mcal{O}_{\mfrak{X}},
\mcal{K}^{\bdot}) \cong
\lim_{n \ar} g^{*} \mcal{K}_{n}^{\bdot} \cong
g^{*} \ul{\Gamma}_{X_{0}} \mcal{K}^{\bdot} \cong
\ul{\Gamma}_{\mrm{disc}} g^{*} \mcal{K}^{\bdot} \]
is a residual complex.
We see that if $\mcal{R}^{\bdot}$ is
any dualizing complex on $X$ then
$\mrm{E} \mrm{R} \ul{\Gamma}_{\mrm{disc}} g^{*} \mcal{R}^{\bdot}$
is dualizing on $\mfrak{X}$.
\end{exa}

We call a formal scheme $\mfrak{X}$ {\em regular} of all its local rings
$\mcal{O}_{\mfrak{X}, x}$ are regular.

\begin{lem} \label{lem5.7}
Suppose $\mfrak{X}$ is a regular formal scheme. Then
$d(x) := \opn{dim} \mcal{O}_{\mfrak{X},x}$
is a bounded codimension function on $\mfrak{X}$.
\end{lem}

\begin{proof}
Let
$\mfrak{U} = \opn{Spf} A \subset \mfrak{X}$ be a connected affine open
set. If $x \in \mfrak{U}$ is the point corresponding to an open prime
ideal $\mfrak{p}$, then
$\widehat{A}_{\mfrak{p}} \cong \widehat{\mcal{O}}_{\mfrak{X}, x}$.
Therefore $A_{\mfrak{p}}$ is a regular local ring.
Now in the adic noetherian ring $A$ any maximal ideal $\mfrak{m}$ is open.
Hence, by \cite{Ma} \S 18 Lemma 5 (III), $A$ is a regular ring,
of finite global dimension equal to its Krull dimension.

Now let
$U := \opn{Spec} A$, so as a topological space, $\mfrak{U} \subset U$
is the closed set defined by any defining ideal $I \subset A$.
Since $U$ is a regular scheme, $\mcal{O}_{U}$ is a dualizing complex
on it. The codimension function $d'$ corresponding to $\mcal{O}_{U}$
satisfies $d'(y) = \opn{dim} \mcal{O}_{U, y}$. Thus
$0 \leq d'(y) \leq \opn{dim} U$. But clearly
$d|_{\mfrak{U}} = d'|_{\mfrak{U}}$. By covering $\mfrak{X}$ with
finitely many such $\mfrak{U}$ this implies
that $d$ is a bounded codimension function.
\end{proof}

\begin{thm} \label{thm5.3}
Suppose $\mfrak{X}$ is a regular formal scheme. Then
$\mrm{R} \ul{\Gamma}_{\mrm{disc}} \mcal{O}_{\mfrak{X}}$
is a dualizing complex on $\mfrak{X}$.
\end{thm}

\begin{proof}
By the proof of Theorem \ref{thm4.2} and known properties of regular
local rings, for any $x \in \mfrak{X}$
\[ \mrm{H}^{q}_{x} \mrm{R} \ul{\Gamma}_{\mrm{disc}}
\mcal{O}_{\mfrak{X}} \cong
\mrm{H}^{q}_{\mfrak{m}_{x}} \widehat{\mcal{O}}_{\mfrak{X}, x} \cong
\begin{cases}
J(x) & \text{ if } q = d(x)\\
0 & \text{ otherwise}
\end{cases} \]
where $\mfrak{m}_{x} \subset \widehat{\mcal{O}}_{\mfrak{X}, x}$
is the maximal ideal, and $J(x)$ is an injective hull of $k(x)$.
Since $d$ is bounded it follows that
$\mcal{K}^{\bdot} :=
\mrm{E} \mrm{R} \ul{\Gamma}_{\mrm{disc}} \mcal{O}_{\mfrak{X}}$
is a bounded complex of injectives in $\msf{QCo}_{\mrm{disc}}(\mfrak{X})$.
Like in the proof of Proposition \ref{prop5.7},
$\mrm{R} \ul{\Gamma}_{\mrm{disc}} \mcal{O}_{\mfrak{X}} \cong
\mcal{K}^{\bdot}$
in $\msf{D}(\msf{Mod}(\mfrak{X}))$.

To complete the proof it suffices to show that for any affine open set
$\mfrak{U} = \opn{Spf} A \subset \mfrak{X}$ the complex
$\mcal{K}^{\bdot}|_{\mfrak{U}}$ is residual on $\mfrak{U}$.
Let $U := \opn{Spec} A$ and let $g : \mfrak{U} \ar U$ be the canonical
morphism Let $U_{0} \subset U$ be the closed set $g(\mfrak{U})$,
so that $\mfrak{U} \cong U_{/ U_{0}}$.
Define
$\mcal{K}^{\bdot}_{U} := \mrm{E} \mcal{O}_{U}$,
which is a residual complex on $U$. Then according to
Proposition \ref{prop4.2}
\[ \mrm{R} \ul{\Gamma}_{\mrm{disc}} \mcal{O}_{\mfrak{U}} \cong
g^{*} \mrm{R} \ul{\Gamma}_{U_{0}} \mcal{O}_{U} \cong
g^{*} \ul{\Gamma}_{U_{0}} \mcal{K}^{\bdot}_{U} . \]
As in Example \ref{exa5.1} this is a dualizing complex, so
$\mcal{K}^{\bdot}|_{\mfrak{U}} \cong
\mrm{E} \mrm{R} \ul{\Gamma}_{\mrm{disc}} \mcal{O}_{\mfrak{U}}$
is a residual complex.
\end{proof}

\begin{rem}
According to \cite{RD} Thm.\ VI.3.1,
if $f : X \ar Y$ is a finite type morphism
between finite dimensional noetherian schemes, and if $\mcal{K}^{\bdot}$
is a residual complex on $Y$, then there is a residual complex
$f^{\triangle} \mcal{K}^{\bdot}$ on $X$. Now suppose
$f : \mfrak{X} \ar \mfrak{Y}$ is a f.f.t.\ morphism and
$f_{n} : X_{n} \ar Y_{n}$ are like in Lemma \ref{lem1.12}.
In the same fashion as in Prop.\ \ref{prop5.8} we set
$f^{\triangle} \mcal{K}^{\bdot} := \lim_{n \ar} f_{n}^{\triangle}
\mcal{K}^{\bdot}_{n}$. This is a residual complex on $\mfrak{X}$.
If $f$ is formally proper then
$\opn{Tr}_{f} = \lim_{n \ar} \opn{Tr}_{f_{n}}$ induces a duality
\[ \mrm{R} f_{*} \mcal{M}^{\bdot} \ar
\mrm{R} \mcal{H}om_{\mfrak{Y}}(\mrm{R} f_{*}
\mrm{R} \mcal{H}om_{\mfrak{X}}(\mcal{M}^{\bdot},
f^{\triangle} \mcal{K}^{\bdot}), \mcal{K}^{\bdot}) \]
for every
$\mcal{M}^{\bdot} \in \msf{D}^{\mrm{b}}(\msf{Coh}(\mfrak{X}))$.
The proofs are standard, given the results of this section.
\end{rem}

% ** section 6 **

\section{Construction of the Complex $\mcal{K}^{\bdot}_{X/S}$}

In this section we work over a regular noetherian base scheme $S$.
We construct the relative residue complex
$\mcal{K}^{\bdot}_{X/S}$ on any finite type $S$-scheme $X$.
The construction is explicit and does not rely on \cite{RD}.

Let $A, B$ be complete local rings, with maximal ideals
$\mfrak{m}, \mfrak{n}$.
Recall that a local homomorphism $\phi : A \ar B$ is called residually
finitely generated if the field extension
$A / \mfrak{m} \ar B / \mfrak{n}$ is finitely generated.
Denote by $\msf{Mod}_{\mrm{disc}}(A)$ the category of $\mfrak{m}$-torsion
$A$-modules (equivalently, modules with $0$-dimensional support).

Suppose
$A \sqbr{\ul{t}} = A \sqbr{t_{1}, \ldots, t_{n}}$
is a polynomial algebra and
$\mfrak{p} \subset A \sqbr{\ul{t}}$ is some maximal ideal.
Then
$A \ar B = \widehat{A \sqbr{\ul{t}}}_{\mfrak{p}}$
is formally smooth of relative dimension $n$ and residually finite.
Let $b_{i} \in B / \mfrak{n}$ be the image of $t_{i}$ and
$\bar{q}_{i} \in
(A / \mfrak{m}) \sqbr{b_{1}, \ldots, b_{i-1}}\sqbr{t_{i}}$
the monic irreducible polynomial of $b_{i}$, of degree $d_{i}$.
Choose a monic lifting
$q_{i} \in A \sqbr{t_{1}, \ldots, t_{i}}$.
Then for a discrete $A$-module $M$ one has
\[ \mrm{H}^{n}_{\mfrak{p}}
\left( \widehat{\Omega}^{n}_{B / A} \otimes_{A} M \right)
\cong \bigoplus_{1 \leq i_{l}} \ \bigoplus_{0 \leq j_{l} < d_{l}}
\gfrac{ t_{1}^{j_{1}} \cdots t_{n}^{j_{n}}
\mrm{d} t_{1} \cdots \mrm{d} t_{n} }
{ q_{1}^{i_{1}} \cdots q_{n}^{i_{n}} } \otimes M . \]
As in \cite{Hg1} Section 7 define the Tate residue
\begin{equation} \label{eqn6.6}
\opn{res}_{t_{1}, \ldots, t_{n}; A, B} :
\mrm{H}^{n}_{\mfrak{p}}
\left( \widehat{\Omega}^{n}_{B / A} \otimes_{A} M \right) \ar M
\end{equation}
by the rule
\[ \gfrac{ t_{1}^{j_{1}} \cdots t_{n}^{j_{n}}
\mrm{d} t_{1} \cdots \mrm{d} t_{n} }
{ q_{1}^{i_{1}} \cdots q_{n}^{i_{n}} } \otimes m \mapsto
\begin{cases}
m & \text{ if } i_{l} = 1, j_{l} = d_{l} - 1 \\
0 & \text{ otherwise}
\end{cases} \]
(cf.\ \cite{Ta}).
Observe that any residually finite homomorphism $A \ar C$ factors into some
$A \ar B = \widehat{A \sqbr{\ul{t}}}_{\mfrak{p}} \ar C$.

\begin{thm} \label{thm6.1} \textup{(Huang)}\
Consider the category $\msf{Loc}$ of complete noetherian local rings and
residually finitely generated local homomorphisms. Then:
\begin{enumerate}
\item For any morphism
$\phi : A \ar B$ in $\msf{Loc}$ there is a functor
\[ \phi_{\#} : \msf{Mod}_{\mrm{disc}}(A)
\ar \msf{Mod}_{\mrm{disc}}(B) . \]
For composable morphisms
$A \xrightarrow{\phi} B \xrightarrow{\psi} C$
there is an isomorphism
$(\psi \phi)_{\#} \cong \psi_{\#} \phi_{\#}$,
and
$(1_{A})_{\#} \cong 1_{\msf{Mod}_{\mrm{disc}}(A)}$.
These data form a pseudofunctor on $\msf{Loc}$ \textup{(}cf.\ \cite{Hg1}
Def.\ \textup{4.1)}.
\item If $\phi : A \ar B$ is formally smooth of relative dimension $q$,
and
$n = \opn{rank} \widehat{\Omega}^{1}_{B / A}$,
then there is an isomorphism, functorial in
$M \in \msf{Mod}_{\mrm{disc}}(A)$,
\[ \phi_{\#} M \cong \mrm{H}^{q}_{\mfrak{n}}(
\widehat{\Omega}^{n}_{B/A} \otimes_{A} M) . \]
\item If $\phi : A \ar B$ is residually finite then
there is an $A$-linear homomorphism, functorial in
$M \in \msf{Mod}_{\mrm{disc}}(A)$,
\[ \opn{Tr}_{\phi} : \phi_{\#} M \ar M , \]
which induces an isomorphism
$\phi_{\#} M \cong \mrm{Hom}^{\mrm{cont}}_{A}(B, M)$.
For composable homomorphisms
$A \xrightarrow{\phi} B \xrightarrow{\psi} C$
one has
$\opn{Tr}_{\psi \phi} = \opn{Tr}_{\phi} \opn{Tr}_{\psi}$
under the isomorphism of part \textup{1}.
\item If $B = \widehat{A \sqbr{\ul{t}}}_{\mfrak{p}}$
then
$\opn{Tr}_{\phi} = \opn{res}_{t_{1}, \ldots, t_{n}; A, B}$
under the isomorphism of part \textup{2}.
\end{enumerate}
\end{thm}

\begin{proof}
Parts 1 and 2 are \cite{Hg1} Thm.\ 6.12. Parts 3 and 4
follow from \cite{Hg1} Section 7.
\end{proof}

\begin{dfn} \label{dfn6.3}
Suppose $L$ is a regular local ring of dimension $q$, with maximal
ideal $\mfrak{r}$. Given a homomorphism $\phi : L \ar A$
in $\msf{Loc}$, define
\[ \mcal{K}(A / L) := \phi_{\#} \mrm{H}^{q}_{\mfrak{r}} L , \]
the {\em dual module of $A$ relative to $L$}.
\end{dfn}

Since $\mrm{H}^{q}_{\mfrak{r}} L$ is an injective hull of the field
$L / \mfrak{r}$, it follows that $\mcal{K}(A / L)$ is an injective hull
of $A / \mfrak{m}$ (cf.\ \cite{Hg1} Corollary 3.10).

\begin{cor} \label{cor6.1}
If $\psi : A \ar B$ is a residually finite homomorphism, then there
is an $A$-linear homomorphism
\[ \opn{Tr}_{\psi} = \opn{Tr}_{B / A} : \mcal{K}(B / L) \ar
\mcal{K}(A / L) . \]
Given another such homomorphism $B \ar C$, one has
$ \opn{Tr}_{C / A} = \opn{Tr}_{B / A} \opn{Tr}_{C / B}$.
\end{cor}

\begin{rem} \label{rem6.1}
One can show that when $L$ is a perfect field,
there is a functorial isomorphism between
$\mcal{K}(A / L) = \phi_{\#} L$ above and the dual module
$\mcal{K}(A)$ of \cite{Ye2}, which was defined via Beilinson completion
algebras.
\end{rem}

Suppose $\pi : \mfrak{X} \ar S$ is a formally finite type (f.f.t.)
formally smooth morphism. According to Proposition \ref{prop1.4},
$\mfrak{X}$ is a regular formal scheme. When we write
$n = \opn{rank} \widehat{\Omega}^{1}_{\mfrak{X}/S}$
we mean that $n$ is a locally constant function
$n : \mfrak{X} \ar \mbb{N}$.

\begin{lem} \label{lem6.1}
Given a f.f.t.\ morphism $\pi : \mfrak{X} \ar S$ and a point
$x \in \mfrak{X}$, let $s := \pi(x)$, and define
\[ d_{S}(x) := \opn{dim} \widehat{\mcal{O}}_{S, s} -
\opn{tr.deg}_{k(s)} k(x) . \]
Then:
\begin{enumerate}
\item $d_{S}$ is a codimension function.
\item If $\pi$ is formally smooth then
\[ d_{S}(x) = \opn{dim} \widehat{\mcal{O}}_{\mfrak{X}, x} -
\opn{rank} \widehat{\Omega}^{1}_{\mfrak{X}/S} . \]
\end{enumerate}
\end{lem}

\begin{proof}
We shall prove 2 first.
Let $L := \widehat{\mcal{O}}_{S, s}$
and
$A := \widehat{\mcal{O}}_{\mfrak{X}, x}$.
By Prop.\ \ref{prop1.4},
\[ \opn{rank} \widehat{\Omega}^{1}_{A/L} =
\opn{dim} A - \opn{dim} L + \opn{tr.deg}_{L / \mfrak{r}} A / \mfrak{m} . \]
We see that $d_{S}$ is the codimension function associated with the
dualizing complex
$\mrm{R} \ul{\Gamma}_{\mrm{disc}} \mcal{O}_{\mfrak{X}}[n]$
(see Theorem \ref{thm5.3}).

As for 1, the property of being a codimension function is local.
But locally there is always a closed immersion
$\mfrak{X} \subset \mfrak{Y}$ with $\mfrak{Y} \ar S$ formally smooth.
\end{proof}

We shall use the codimension function $d_{S}$ by default.

\begin{dfn} \label{dfn6.1}
Let $\pi : \mfrak{X} \ar S$ be a formally finite type morphism. Given a
point $x \in \mfrak{X}$, consider
$\phi : L = \widehat{\mcal{O}}_{S ,\pi(x)} \ar
A = \widehat{\mcal{O}}_{\mfrak{X}, x}$,
which is a morphism in $\msf{Loc}$.
Since $L$ is a regular local ring, the dual module $\mcal{K}(A / L)$
is defined.
Let $\mcal{K}_{\mfrak{X} / S}(x)$ be the quasi-coherent sheaf
which is constant on
$\overline{\{x\}}$ with group of sections $\mcal{K}(A / L)$,
and define
\[ \mcal{K}_{\mfrak{X} / S}^{q} :=
\bigoplus_{d_{S}(x) = q} \mcal{K}_{\mfrak{X} / S}(x) . \]
\end{dfn}

In Theorem \ref{thm6.2} we are going to prove that on
the graded sheaf $\mcal{K}_{X / S}^{\bdot}$
there is a canonical coboundary operator $\delta$ which makes it into
residual complex.

\begin{dfn} \label{dfn6.2}
Let $f : \mfrak{X} \ar \mfrak{Y}$ be a morphism of formal schemes over $S$.
Define a homomorphism of graded $\mcal{O}_{\mfrak{Y}}$-modules
$\opn{Tr}_{f} : f_{*} \mcal{K}_{\mfrak{X} / S}^{\bdot}  \ar
\mcal{K}_{\mfrak{Y} / S}^{\bdot}$
as follows. If $x \in \mfrak{X}$ is closed in its fiber and $y = f(x)$,
then
$A = \widehat{\mcal{O}}_{\mfrak{Y}, y} \ar
B = \widehat{\mcal{O}}_{\mfrak{X}, x}$
is a residually finite $L$-algebra homomorphism.
The homomorphism
$\opn{Tr}_{B / A} : \mcal{K}(B / L) \ar \mcal{K}(A / L)$
of Cor.\ \ref{cor6.1} gives a map of sheaves
\[ \opn{Tr}_{f} : f_{*} \mcal{K}_{\mfrak{X} / S}(x)  \ar
\mcal{K}_{\mfrak{Y} / S}(y) . \]
If $x$ is not closed in its fiber, we let $\opn{Tr}_{f}$ vanish on
$f_{*} \mcal{K}_{\mfrak{X} / S}(x)$.
\end{dfn}

\begin{prop} \label{prop6.3}
\begin{enumerate}
\item $\opn{Tr}_{f}$ is functorial: if $g : \mfrak{Y} \ar \mfrak{Z}$
is another morphism, then
$\opn{Tr}_{gf} = \opn{Tr}_{g} \opn{Tr}_{f}$.
\item If $f$ is formally finite \textup{(}see Def.\
\textup{\ref{dfn1.1})},
then $\opn{Tr}_{f}$ induces an isomorphism
of graded sheaves
\[ f_{*} \mcal{K}_{\mfrak{X} / S}^{\bdot} \cong
\mcal{H}om_{\mfrak{Y}}(f_{*} \mcal{O}_{\mfrak{X}},
\mcal{K}_{\mfrak{Y} / S}^{\bdot}) . \]
\item If $g : \mfrak{U} \ar \mfrak{X}$ is an open immersion, then
there is a natural isomorphism
$\mcal{K}_{\mfrak{U} / S}^{\bdot} \cong
g^{*} \mcal{K}_{\mfrak{X} / S}^{\bdot}$.
\end{enumerate}
\end{prop}

\begin{proof}
Part 3 is trivial. Part 1 is a consequence of Cor.\ \ref{cor6.1}.
As for part 2, $f$ is an affine morphism, and fibers of $f$ are all
finite, so all points of $X$ are closed in their fibers.
\end{proof}

Suppose $\ul{a} = (a_{1}, \ldots, a_{n})$ is a sequence of elements
in the noetherian ring $A$.
Let us write $\tilde{\mbf{K}}^{\bdot}_{\infty}(\ul{a})$
for the subcomplex $\mbf{K}^{\geq 1}_{\infty}(\ul{a})$, so we get an
exact sequence
\begin{equation}
0 \ar \tilde{\mbf{K}}^{\bdot}_{\infty}(\ul{a}) \ar
\mbf{K}^{\bdot}_{\infty}(\ul{a}) \ar A \ar  0 .
\end{equation}
For any
$M^{\bdot} \in \msf{D}^{+}(\msf{Mod}(A))$
let $\mcal{M}^{\bdot}$ be the complex of
sheaves $\mcal{O}_{X} \otimes M^{\bdot}$ on
$X := \opn{Spec} A$, and let $U \subset X$ be the open set
$\bigcup \{ a_{i} \neq 0 \}$. Then
\[ \mrm{R} \Gamma(U, \mcal{M}^{\bdot}) \cong
\tilde{\mbf{K}}^{\bdot}_{\infty}(\ul{a})[1] \otimes M^{\bdot} \]
in $\msf{D}(\msf{Mod}(A))$. In fact
$\tilde{\mbf{K}}^{\bdot}_{\infty}(\ul{a}) \otimes \mcal{O}_{X}$
is a shift by $1$ of the \v{C}ech complex corresponding to the open
cover of $U$.

\begin{lem} \label{lem6.2}
Let $A$ be an adic noetherian ring and
$M^{\bdot} \in \msf{D}^{+}(\msf{Mod}(A))$.
Define $\mfrak{U} := \opn{Spf} A$ and
$\mcal{M}^{\bdot} := \mcal{O}_{\mfrak{U}} \otimes M^{\bdot}$.
\begin{enumerate}
\item Let $x \in \mfrak{U}$ with corresponding open prime ideal
$\mfrak{p} \subset A$. Suppose the sequence
$\ul{a}$ generates $\mfrak{p}$. Then
\[ \mrm{R} \Gamma_{x} \mrm{R} \ul{\Gamma}_{\mrm{disc}}
\mcal{M}^{\bdot} \cong
\mrm{R} \Gamma_{\mfrak{p}} M^{\bdot}_{\mfrak{p}} \cong
\mbf{K}^{\bdot}_{\infty}(\ul{a}) \otimes M^{\bdot}_{\mfrak{p}} \]
in $\msf{D}^{+}(\msf{Mod}(A_{\mfrak{p}})).$
\item Suppose $y \in \mfrak{U}$ is an immediate specialization of $x$,
and its ideal $\mfrak{q}$ has generators $\ul{a}, \ul{b}$. Then
\[ \mrm{R} \Gamma_{x} \mrm{R} \ul{\Gamma}_{\mrm{disc}}
\mcal{M}^{\bdot} \cong
\mbf{K}^{\bdot}_{\infty}(\ul{a}) \otimes
\tilde{\mbf{K}}^{\bdot}_{\infty}(\ul{b})[1] \otimes
M^{\bdot}_{\mfrak{q}} \]
in $\msf{D}^{+}(\msf{Mod}(A_{\mfrak{q}})).$
\item Assume $d$ is a codimension function on $\mfrak{U}$.
Then in the Cousin complex
$\mrm{E} \mrm{R} \ul{\Gamma}_{\mrm{disc}} \mcal{M}^{\bdot}$
the map
\[ \mrm{H}^{d(x)}_{x} \mrm{R} \ul{\Gamma}_{\mrm{disc}} \mcal{M}^{\bdot}
\ar \mrm{H}^{d(y)}_{y} \mrm{R} \ul{\Gamma}_{\mrm{disc}} \mcal{M}^{\bdot} \]
is given by applying $\mrm{H}^{d(y)}$ to
\[ \left( \mbf{K}^{\bdot}_{\infty}(\ul{a}) \otimes
\tilde{\mbf{K}}^{\bdot}_{\infty}(\ul{b}) \ar
\mbf{K}^{\bdot}_{\infty}(\ul{a}, \ul{b}) \right) \otimes
M^{\bdot}_{\mfrak{q}} . \]
\end{enumerate}
\end{lem}

\begin{proof}
Part 1 follows immediately from formula (\ref{eqn4.1}).
Parts 2 and 3 are true because
$\opn{Spec} (A / \mfrak{p})_{\mfrak{q}} =
\{ \mfrak{p}, \mfrak{q} \}$.
\end{proof}

As a warm up for Thm.\ \ref{thm6.2}, here is:

\begin{prop} \label{prop6.4}
If $\pi : \mfrak{X} \ar S$ is formally smooth, with
$n = \opn{rank} \widehat{\Omega}^{1}_{\mfrak{X} / S}$,
then there is a canonical isomorphism of graded sheaves
\[ \mcal{K}_{\mfrak{X} / S}^{\bdot} \cong
\mrm{E} \mrm{R} \ul{\Gamma}_{\mrm{disc}}
\widehat{\Omega}^{n}_{\mfrak{X} / S}[n] . \]
This makes $\mcal{K}_{\mfrak{X} / S}^{\bdot}$ into a residual complex.
\end{prop}

\begin{proof}
Take any point $x$, and with the notation of Def.\
\ref{dfn6.1} let $p := \opn{dim} L$ and $q := \opn{dim} A$.
Then by Lemma \ref{lem6.2} part 1 and \cite{Hg1} Proposition 2.6
we have a canonical isomorphism
\[ \mrm{H}^{d(x)}_{x}
\mrm{R} \ul{\Gamma}_{\mrm{disc}} \widehat{\Omega}^{n}_{\mfrak{X} / S}[n]
\cong
\mrm{H}^{q}_{\mfrak{m}} \widehat{\Omega}^{n}_{A / L} \cong
\mrm{H}^{q-p}_{\mfrak{m}} \left( \widehat{\Omega}^{n}_{A / L}
\otimes_{L} \mrm{H}^{p}_{\mfrak{r}} L \right)
\cong \mcal{K}(A / L) . \]
According to Theorem \ref{thm5.3} and Proposition \ref{prop5.7},
$\mrm{E} \mrm{R} \ul{\Gamma}_{\mrm{disc}}
\widehat{\Omega}^{n}_{\mfrak{X} / S}[n]$
is a residual complex.
\end{proof}

In particular taking $\mfrak{X} = S$ we get
$\mcal{K}_{S / S}^{\bdot} = \mrm{E} \mcal{O}_{S}$.

\begin{lem} \label{lem6.3}
Suppose $X \subset \mfrak{X}$ and $X \subset \mfrak{Y}$ are
s.f.e.'s and $f : \mfrak{X} \ar \mfrak{Y}$ is a morphism of
embeddings.
Then
$\opn{Tr}_{f} : \mcal{K}^{\bdot}_{\mfrak{X}} \ar
\mcal{K}^{\bdot}_{\mfrak{Y}}$
is a homomorphism of complexes.
\end{lem}

\begin{proof}
Factoring $f$ through $(\mfrak{X} \times_{S} \mfrak{Y})_{/ X}$
we can assume that $f$ is either a closed immersion, or that it is
 formally smooth. At any rate $f$ is an affine morphism, so we can take
$\mfrak{X} = \opn{Spf} B$, $\mfrak{Y} = \opn{Spf} A$ and
$S = \opn{Spec} L$.
By Theorem  \ref{thm2.2} we can suppose one of the following holds:
(i) $B \cong A [\sqbr{\ul{t}}]$ for a sequence of indeterminates
$\ul{t} = (t_{1}, \ldots, t_{l})$, and $A \ar B$ is the inclusion;
or
(ii) $A \cong B [\sqbr{\ul{t}}]$ and $A \ar B$ is the projection modulo
$\ul{t}$. We shall treat each case separately.

\noindent (i)\
Choose generators $\ul{a}$ for a defining ideal of $A$. Let
$m := \opn{rank} \widehat{\Omega}^{1}_{A / L}$
and
$n := \opn{rank} \widehat{\Omega}^{1}_{B / L}$,
so $n = m + l$. Define an $A$-linear map
$\rho : \mbf{K}^{\bdot}_{\infty}(\ul{t}) \otimes
\widehat{\Omega}^{l}_{B / A}[l] \ar A$
by
$\rho(\ul{t}^{(-1, \ldots, -1)} \mrm{d} \ul{t}) = 1$
and
$\rho(\ul{t}^{\ul{i}}\, \mrm{d} \ul{t}) = 0$ if
$\ul{i} \neq (-1, \ldots, -1)$.
Extend $\rho$ linearly to
\[ \rho : \mbf{K}^{\bdot}_{\infty}(\ul{a}, \ul{t}) \otimes
\widehat{\Omega}^{n}_{B / L}[n] \ar
\mbf{K}^{\bdot}_{\infty}(\ul{a}) \otimes
\widehat{\Omega}^{m}_{A / L}[m] . \]
This $\rho$ sheafifies to give a map of complexes in $\msf{Ab}(X)$
\[ \tilde{\rho} :
\mbf{K}^{\bdot}_{\infty}(\ul{a}, \ul{t}) \otimes
\widehat{\Omega}^{n}_{\mfrak{X} / S}[n] \ar
\mbf{K}^{\bdot}_{\infty}(\ul{a}) \otimes
\widehat{\Omega}^{m}_{\mfrak{Y} / S}[m] . \]
By Lemma \ref{lem6.2} and \cite{Hg1} \S 5, for any point $x \in X$,
$\mrm{H}^{d(x)}_{x}(\tilde{\rho})$ recovers
$\opn{Tr}_{f} : \mcal{K}_{\mfrak{X} / S}(x) \ar
\mcal{K}_{\mfrak{Y} / S}(x)$.
Thus $\opn{Tr}_{f} = \mrm{E}(\tilde{\rho})$ is a homomorphism of complexes.

\noindent (ii)\
Now $l = m - n$. Take $\ul{a}$ to be generators of a defining ideal of
$B$. Define a $B$-linear map
$\rho' : B \ar \mbf{K}^{\bdot}_{\infty}(\ul{t}) \otimes
\widehat{\Omega}^{l}_{A / L}[l]$
by
$\rho'(1) = \ul{t}^{(-1, \ldots, -1)} \mrm{d} \ul{t}$.
Extend $\rho$ linearly to
\[ \rho' : \mbf{K}^{\bdot}_{\infty}(\ul{a}) \otimes
\widehat{\Omega}^{n}_{B / L}[n] \ar
\mbf{K}^{\bdot}_{\infty}(\ul{a}, \ul{t}) \otimes
\widehat{\Omega}^{m}_{A / L}[m] . \]
Again this extends to a map of complexes of sheaves
$\tilde{\rho}'$ in $\msf{Ab}(X)$,
and checking punctually we see that
$\opn{Tr}_{f} = \mrm{E}(\tilde{\rho}')$.
\end{proof}

\begin{thm} \label{thm6.2}
Suppose $X \ar S$ is a finite type morphism.
There is a unique operator
$\delta : \mcal{K}_{X / S}^{q} \ar \mcal{K}_{X / S}^{q+1}$,
satisfying the following local condition:\\[2mm]
\textup{\bf (LE)}\ \blnk{4mm}
\begin{minipage}{11cm}
Suppose $U \subset X$ is
an open subset, and $U \subset \mfrak{U}$ is a smooth formal embedding.
By Proposition \textup{\ref{prop6.3}} there is an inclusion of graded
$\mcal{O}_{U}$-modules
$\mcal{K}_{X / S}^{\bdot}|_{U} \subset
\mcal{K}_{\mfrak{U} / S}^{\bdot}$.
Then $\delta|_{U}$ is compatible with the coboundary operator on
$\mcal{K}_{\mfrak{U} / S}^{\bdot}$ coming from Proposition
\textup{\ref{prop6.4}}.
\end{minipage}\\[2mm]
Moreover $(\mcal{K}_{X / S}^{\bdot}, \delta)$
is a residual complex on $X$.
\end{thm}

\begin{proof}
Define $\delta|_{U}$ using {\bf LE}. According to Lemma \ref{lem6.3},
$\delta|_{U}$ is independent of $\mfrak{U}$, so it glues.
We get a bounded complex of quasi-coherent injectives on $X$.
By Proposition \ref{prop6.4} it follows that it is residual.
\end{proof}

\begin{rem}
This construction of $\mcal{K}_{X / S}^{\bdot}$ actually allows
a computation of the operator $\delta$, given the data of a local
embedding. The formula is in part 3 of Lemma \ref{lem6.2},
with $M^{\bdot} = \widehat{\Omega}^{n}_{A / L}[n]$.
The formula for changing the embedding can be extracted from the
proof of Lemma \ref{lem6.3}. Of course when
$\opn{rank} \widehat{\Omega}^{1}_{\mfrak{X} / S}$ is high
these computations can be nasty.
\end{rem}

\begin{rem} \label{rem6.6}
The recent papers \cite{Hg2}, \cite{Hg3} and \cite{LS2} also
use the local theory of \cite{Hg1} as a starting point for explicit
constructions of Grothendieck Duality. Their constructions are more
general than ours: Huang constructs $f^{!} \mcal{M}^{\bdot}$ for a finite
type morphism $f : X \ar Y$ and a residual complex complex
$\mcal{M}^{\bdot}$; and Lipman-Sastry even allow $\mcal{M}^{\bdot}$
to be any Cousin complex.
\end{rem}

% ** section 7 **

\section{The Trace for Finite Morphisms}

In this section we prove that
$\opn{Tr}_{f}$ is a homomorphism of complexes when $f$ is a finite
morphism. The proof is by a self contained calculation involving
Koszul complexes and a comparison of global and local Tate residue maps.
In Theorem \ref{thm7.4} we compare the complex
$\mcal{K}_{X / S}^{\bdot}$ to the
sheaf of regular differentials of Kunz-Waldi.
Throughout $S$ is a regular noetherian scheme.

\begin{thm} \label{thm7.6}
Suppose $f : X \ar Y$ is finite. Then
$\opn{Tr}_{f} : f_{*} \mcal{K}_{X / S}^{\bdot} \ar
 \mcal{K}_{Y / S}^{\bdot}$
is a homomorphism of complexes.
\end{thm}

The proof appears after some preparatory work, based on and
inspired by \cite{Hg1} \S7.

\begin{rem} \label{rem7.1}
In Section 8 we prove a much stronger result, namely Corollay \ref{cor8.1},
but its proof is indirect and relies on the Residue Theorem of \cite{RD}
Chapter VII. We have decided to include Theorem \ref{thm7.6} because of
its direct algebraic proof.
\end{rem}

Let $A$ be an adic noetherian ring with defining ideal $\mfrak{a}$.
Suppose $p \in A \sqbr{t}$ is a monic polynomial of degree $e > 0$. Define
an $A$-algebra
\begin{equation} \label{eqn7.3}
B := \lim_{\leftarrow i} A \sqbr{t} / A \sqbr{t} \cdot p^{i} .
\end{equation}
Let
$\mfrak{b} :=B  \mfrak{a} + B p$; then
$B \cong \lim_{\leftarrow i} B / \mfrak{b}^{i}$, so that $B$ is an adic
ring with the $\mfrak{b}$-adic topology.
The homomorphism $\phi : A \ar B$ is f.f.t.\ and formally smooth, and
$\widehat{\Omega}^{1}_{B / A} = B \cdot \mrm{d} t$.
Furthermore $p \in B$ is a non-zero-divisor,
and by long division we obtain an isomorphism
\begin{equation} \label{eqn7.2}
\mrm{H}^{1}_{(p)} B =
\mrm{H}^{1} \left(\mbf{K}^{\bdot}_{\infty}(p) \otimes B \right) \cong
\bigoplus_{1 \leq i}\ \bigoplus_{0 \leq j < e} A \cdot
\gfrac{t^{j}}{p^{i}} .
\end{equation}

Define an $A$-linear homomorphism
$\opn{Res}_{B / A} : \mrm{H}^{1}_{(p)} \widehat{\Omega}^{1}_{B / A}
\ar A$
by
\[ \opn{Res}_{B / A} \left( \gfrac{t^{j} \mrm{d}t}{p^{i}} \right) :=
\begin{cases}
1 & \text{ if } i=1, j=e-1 \\
0 & \text{ otherwise} .
\end{cases} \]
We call $\opn{Res}_{B / A}$ the {\em global Tate residue}. It
gives rise to a map of complexes in $\msf{Mod}(A)$:
\begin{equation} \label{eqn7.9}
\opn{Res}_{B / A} : \mbf{K}^{\bdot}_{\infty}(p)[1] \otimes
\widehat{\Omega}^{1}_{B / A} \ar A .
\end{equation}
Note that both the algebra $B$ and the map $\opn{Res}_{B / A}$
depend on $t$ and $p$.

Suppose $\mfrak{q} \subset B$ is an open prime ideal and
$\mfrak{p} = \phi^{-1}(\mfrak{q}) \subset A$. Then the local homomorphism
$\phi_{\mfrak{q}} : \widehat{A}_{\mfrak{p}} \ar
\widehat{B}_{\mfrak{q}}$
is formally smooth of relative dimension $1$ and residually finite.
Let
$\tilde{\mfrak{q}} := \mfrak{q} \cap  \widehat{A}_{\mfrak{p}} \sqbr{t}$,
and denote by $\bar{\mfrak{q}}$ the image of $\tilde{\mfrak{q}}$
in $k(\mfrak{p}) \sqbr{t}$, so
$k(\mfrak{p}) \sqbr{t} / \bar{\mfrak{q}} = k(\mfrak{q})$.
For a polynomial $q \in \widehat{A}_{\mfrak{p}} \sqbr{t}$
let $\bar{q}$ be its image in $k(\mfrak{p}) \sqbr{t}$.
Suppose $q$ satisfies:
\begin{equation} \label{eqn7.5}
q \text{ is monic, and the ideal }
(\bar{q}) \subset k(\mfrak{p}) \sqbr{t} \text{ is }
\bar{\mfrak{q}}\text{-primary.}
\end{equation}
Then
$\widehat{B}_{\mfrak{q}} \cdot \mfrak{q} =
\sqrt{\widehat{B}_{\mfrak{q}} \cdot (\mfrak{p}, q)} \subset
\widehat{B}_{\mfrak{q}}$,
and
\[ \widehat{B}_{\mfrak{q}} \cong
\lim_{\leftarrow i} \widehat{A}_{\mfrak{p}} \sqbr{t} / \tilde{\mfrak{q}}^{i}
\cong \lim_{\leftarrow i}
\widehat{A}_{\mfrak{p}} \sqbr{t} / \widehat{A}_{\mfrak{p}} \sqbr{t}
\cdot q^{i} . \]
Hence $q$ is a non-zero-divisor in $\widehat{B}_{\mfrak{q}}$ and
$\widehat{B}_{\mfrak{q}} / \widehat{B}_{\mfrak{q}} \cdot q$ is a free
$\widehat{A}_{\mfrak{p}}$-module with basis
$1, t, \ldots, t^{d-1}$, where $d = \opn{deg} q$. We see that
a decomposition like (\ref{eqn7.2}) exists for
$\mrm{H}^{1}_{(q)} \widehat{B}_{\mfrak{q}}$.

Suppose we are given a discrete $\widehat{A}_{\mfrak{p}}$-module $M$.
Then one gets
\[ \mrm{H}^{1}_{\mfrak{q}} \left(
\widehat{\Omega}^{1}_{\widehat{B}_{\mfrak{q}} / \widehat{A}_{\mfrak{p}}}
\otimes_{\widehat{A}_{\mfrak{p}}} M \right) \cong
\left( \mrm{H}^{1}_{(q)}
\widehat{\Omega}^{1}_{\widehat{B}_{\mfrak{q}} / \widehat{A}_{\mfrak{p}}}
\right) \otimes_{\widehat{A}_{\mfrak{p}}} M
\cong
\bigoplus_{1 \leq i}\ \bigoplus_{0 \leq j < d}
\gfrac{t^{j} \mrm{d} t}{q^{i}} \otimes M  \]
(cf.\ \cite{Hg1} pp.\ 41-42). Define the {\em local Tate residue map}
\[ \opn{Res}_{\widehat{B}_{\mfrak{q}} / \widehat{A}_{\mfrak{p}}} :
\mrm{H}^{1}_{\mfrak{q}} \left(
\widehat{\Omega}^{1}_{\widehat{B}_{\mfrak{q}} / \widehat{A}_{\mfrak{p}}}
\otimes_{\widehat{A}_{\mfrak{p}}} M \right) \ar M \]
by
\[ \opn{Res}_{\widehat{B}_{\mfrak{q}} / \widehat{A}_{\mfrak{p}}}
\left( \gfrac{t^{j} \mrm{d}t \otimes m}{q^{i}}
\right) :=
\begin{cases}
m & \text{ if } i=1, j=d-1 \\
0 & \text{ otherwise} .
\end{cases} \]
Clearly $\opn{Res}_{\widehat{B}_{\mfrak{q}} / \widehat{A}_{\mfrak{p}}}$
is functorial in $M$, and it depends on $t$.

\begin{lem} \label{lem7.7}
$\opn{Res}_{\widehat{B}_{\mfrak{q}} / \widehat{A}_{\mfrak{p}}}$
is independent of $q$. It coincides with the residue map
$\opn{res}_{t; \widehat{B}_{\mfrak{q}} / \widehat{A}_{\mfrak{p}}}$
of \textup{(\ref{eqn6.6})},
i.e.\ of \cite{Hg1} Definition \textup{8.1}.
\end{lem}

\begin{proof}
Suppose the polynomials
$q_{1}, q_{2} \in \widehat{A}_{\mfrak{p}} \sqbr{t}$
satisfy (\ref{eqn7.5}). Then so does $q_{3} := q_{1} q_{2}$.
Let $\opn{deg} q_{h} = d_{h}$, and let
$\opn{Res}_{\widehat{B}_{\mfrak{q}} / \widehat{A}_{\mfrak{p}}; q_{h}}$
be the residue map determined by $q_{h}$.
Pick any $1 \leq i$ and $0 \leq j < d_{1}$, and write
$q_{2}^{i} = \sum_{l = 0}^{i d_{2}} a_{l} t^{l}$, so $a_{i d_{2}} = 1$.
By the rules for manipulating generalized fractions (cf.\ \cite{Hg1}
\S 1) we have
\begin{equation} \label{eqn7.1}
\opn{Res}_{\widehat{B}_{\mfrak{q}} / \widehat{A}_{\mfrak{p}}; q_{3}}
\left( \gfrac{t^{j} \mrm{d} t \otimes m}{q_{1}^{i}} \right) =
\sum_{l = 0}^{i d_{2}}
\opn{Res}_{\widehat{B}_{\mfrak{q}} / \widehat{A}_{\mfrak{p}}; q_{3}}
\left(
\gfrac{t^{l + j} \mrm{d} t \otimes a_{l} m}{q_{3}^{i}} \right) .
\end{equation}
If $i \geq 2$ or $j \leq d_{1} - 2$ one has
$l + j \leq i d_{3} -2$,
and therefore each summand of the right side of (\ref{eqn7.1}) is $0$.
When $i = 1$ and $j = d_{1} - 1$ the only possible nonzero
residue there is for $l = d_{2}$, and this residue is $m$.
We conclude that
$\opn{Res}_{\widehat{B}_{\mfrak{q}} / \widehat{A}_{\mfrak{p}}; q_{3}} =
\opn{Res}_{\widehat{B}_{\mfrak{q}} / \widehat{A}_{\mfrak{p}}; q_{1}}$.
Clearly also
$\opn{Res}_{\widehat{B}_{\mfrak{q}} / \widehat{A}_{\mfrak{p}}; q_{3}} =
\opn{Res}_{\widehat{B}_{\mfrak{q}} / \widehat{A}_{\mfrak{p}}; q_{2}}$.

If we take $q$ such that $(\bar{q}) = \bar{\mfrak{q}}$, this is by
definition the residue map of (\ref{eqn6.6}).
\end{proof}

\begin{lem} \label{lem7.8}
Let $F$ be the set of prime ideals in $B / (p)$ lying over
$\mfrak{p}$. Then for any
$M \in \msf{Mod}_{\mrm{disc}}(\widehat{A}_{\mfrak{p}})$ one has
\[ \left( \mrm{H}^{1}_{(p)} \widehat{\Omega}^{1}_{B / A} \right)
\otimes_{A} M
\cong
\bigoplus_{\mfrak{q}' \in F}
\mrm{H}^{1}_{\mfrak{q}'} \left(
\widehat{\Omega}^{1}_{\widehat{B}_{\mfrak{q}'} / \widehat{A}_{\mfrak{p}}}
\otimes_{\widehat{A}_{\mfrak{p}}} M \right) , \]
and w.r.t.\ this isomorphism,
\[ \opn{Res}_{B / A} \otimes 1 = \sum_{\mfrak{q}' \in F}
\opn{Res}_{\widehat{B}_{\mfrak{q}'} / \widehat{A}_{\mfrak{p}}} . \]
\end{lem}

\begin{proof}
The isomorphism of modules is not hard to see. Let
$\bar{p} = \prod_{\mfrak{q}' \in F} \bar{p}_{\mfrak{q}'}$
be the primary decomposition in $k(\mfrak{p}) \sqbr{t}$
(all the $\bar{p}_{\mfrak{q}'}$ monic). By Hensel's Lemma this
decomposition lifts to
$p = \prod_{\mfrak{q}' \in F} p_{\mfrak{q}'}$
in $\widehat{A}_{\mfrak{p}} \sqbr{t}$.
Since each polynomial $p_{\mfrak{q}'}$ satisfies condition (\ref{eqn7.5})
for the prime ideal $\mfrak{q}'$, we can use it to calculate
$\opn{Res}_{\widehat{B}_{\mfrak{q}'} / \widehat{A}_{\mfrak{p}}}$.
\end{proof}

\begin{proof} (of Thm.\ \ref{thm7.6})\
This claim is local on $Y$, so we may assume $X$, $Y$ and $S$ are affine,
say $X = \opn{Spec} \bar{B}$, $Y = \opn{Spec} \bar{A}$ and
$S = \opn{Spec} L$.
By the functoriality of $\opn{Tr}$ we can assume
$\bar{B} = \bar{A} \sqbr{b}$
for some element $b \in \bar{B}$. It will suffice to find suitable
s.f.e.'s $X \subset \mfrak{X}$ and $Y \subset \mfrak{Y}$
with a morphism
$\widehat{f} : \mfrak{X} \ar \mfrak{Y}$
extending $f$, and to check that
$\opn{Tr}_{\widehat{f}} : \widehat{f}_{*}
\mcal{K}^{\bdot}_{\mfrak{X} / S} \ar \mcal{K}^{\bdot}_{\mfrak{Y} / S}$
commutes with $\delta$.

Pick any s.f.e.\
$Y \subset \mfrak{Y} = \opn{Spf} A$, so
$\mfrak{a} := \opn{Ker}(A \ar \bar{A})$ is a defining ideal.
Let $A \sqbr{t} \ar \bar{B}$ be the homomorphism $t \mapsto b$.
Choose any monic polynomial $p(t) \in A \sqbr{t}$ s.t.\ $p(b) = 0$, and
define the adic ring $B$ as in formula (\ref{eqn7.3}). So
$\mfrak{X} := \opn{Spf} B$ is the s.f.e.\ of $X$ we want.

Let $(y_{0}, y_{1})$ be an immediate specialization pair in $Y$, and let
$F_{i} := f^{-1}(y_{i}) \subset X$.
Let $\mfrak{p}_{0} \subset \mfrak{p}_{1} \subset A$ be the prime ideals
corresponding to $(y_{0}, y_{1})$.
Pick a sequence of generators $\ul{a}$ for $\mfrak{p}_{0}$, and
generators $(\ul{a}, \ul{a}')$ for $\mfrak{p}_{1}$.
Let $m := \opn{rank} \widehat{\Omega}^{1}_{A / L}$.

Consider the commutative diagram of complexes
\[ \begin{CD}
\mbf{K}^{\bdot}_{\infty}(\ul{a}, p)[1] \otimes
\tilde{\mbf{K}}^{\bdot}_{\infty}(\ul{a}') \otimes
(\widehat{\Omega}^{m + 1}_{B / L})_{\mfrak{p}_{1}}
@> \opn{Res}_{B / A} \otimes 1 >>
\mbf{K}^{\bdot}_{\infty}(\ul{a}) \otimes
\tilde{\mbf{K}}^{\bdot}_{\infty}(\ul{a}') \otimes
(\widehat{\Omega}^{m}_{A / L})_{\mfrak{p}_{1}} \\
@VVV @VVV \\
\mbf{K}^{\bdot}_{\infty}(\ul{a}, \ul{a}', p)[1] \otimes
(\widehat{\Omega}^{m + 1}_{B / L})_{\mfrak{p}_{1}}
@> \opn{Res}_{B / A} \otimes 1 >>
\mbf{K}^{\bdot}_{\infty}(\ul{a}, \ul{a}') \otimes
(\widehat{\Omega}^{m}_{A / L})_{\mfrak{p}_{1}}
\end{CD} \]
gotten from tensoring the map $\opn{Res}_{B / A}$ of (\ref{eqn7.9})
with
$A_{\mfrak{p}_{1}} \otimes \widehat{\Omega}^{m}_{A / L}$
and the various $\mbf{K}^{\bdot}_{\infty}$.
Applying $\mrm{H}^{i}$ to this diagram, where
$i := \opn{dim} \widehat{A}_{\mfrak{p}_{1}}$,
and using Lemmas \ref{lem6.2} and \ref{lem7.8} we obtain a
commutative diagram
\[ \begin{CD}
\bigoplus_{\mfrak{q}_{0} \in F_{0}}
\mrm{H}^{1}_{\mfrak{q}_{0}} \left( \widehat{\Omega}^{1}_{
\widehat{B}_{\mfrak{q}_{0}} / \widehat{A}_{\mfrak{p}_{0}}}
\otimes
\mrm{H}^{i - 1}_{\mfrak{p}_{0}} \widehat{\Omega}^{m}_{
\widehat{A}_{\mfrak{p}_{0}} / L} \right)
@> \sum \opn{Res} >>
\mrm{H}^{i - 1}_{\mfrak{p}_{0}} \widehat{\Omega}^{m}_{
\widehat{A}_{\mfrak{p}_{0}} / L} \\
@VVV @VVV \\
\bigoplus_{\mfrak{q}_{1} \in F_{1}}
\mrm{H}^{1}_{\mfrak{q}_{1}} \left( \widehat{\Omega}^{1}_{
\widehat{B}_{\mfrak{q}_{1}} / \widehat{A}_{\mfrak{p}_{1}}}
\otimes
\mrm{H}^{i}_{\mfrak{p}_{1}} \widehat{\Omega}^{m}_{
\widehat{A}_{\mfrak{p}_{1}} / L} \right)
@> \sum \opn{Res} >>
\mrm{H}^{i}_{\mfrak{p}_{1}} \widehat{\Omega}^{m}_{
\widehat{A}_{\mfrak{p}_{1}} / L} .
\end{CD} \]
In this diagram
$\opn{Res} = \opn{Res}_{
\widehat{B}_{\mfrak{q}_{0}} / \widehat{A}_{\mfrak{p}_{0}}}$
etc. Using the definitions this is the same as
\[ \begin{CD}
\bigoplus_{x_{0} \in F_{0}}
f_{*} \mcal{K}_{\mfrak{X} / S}(x_{0})
@> \opn{Tr}_{f} >>
\mcal{K}_{\mfrak{Y} / S}(y_{0}) \\
@V \delta VV @V \delta VV \\
\bigoplus_{x_{1} \in F_{1}}
f_{*} \mcal{K}_{\mfrak{X} / S}(x_{1})
@> \opn{Tr}_{f} >>
\mcal{K}_{\mfrak{Y} / S}(y_{1}) .
\end{CD} \]
\end{proof}

According to \cite{KW}, if $\pi : X \ar S$ is equidimensional of
dimension $n$ and generically smooth, and $X$ is integral,
then the {\em sheaf of regular
differentials} $\tilde{\omega}^{n}_{X/S}$ (relative to the DGA
$\mcal{O}_{S}$) exists. It is a coherent subsheaf of
$\Omega^{n}_{k(X)/k(S)}$.

\begin{thm} \label{thm7.4}
Suppose $\pi : X \ar S$ is equidimensional of
dimension $n$ and generically smooth, and $X$ is integral. Then
$\mcal{K}^{-n}_{X/S} = \Omega^{n}_{k(X)/k(S)}$, and
\[ \tilde{\omega}^{n}_{X/S} = \mrm{H}^{-n} \mcal{K}^{\bdot}_{X/S} . \]
\end{thm}

First we need:

\begin{lem} \label{lem7.1}
Suppose $L_{0} \ar A_{0} \ar B_{0}$ are finitely generated field
extensions, with
$L_{0} \ar A_{0}$ and $L_{0} \ar B_{0}$ separable, $A_{0} \ar B_{0}$
finite, and
$\opn{tr.deg}_{L_{0}} A_{0} = n$. Then
$\mcal{K}(A_{0} / L_{0}) = \Omega^{n}_{A_{0} / L_{0}}$,
$\mcal{K}(B_{0} / L_{0}) = \Omega^{n}_{B_{0} / L_{0}}$,
and
$\opn{Tr}_{B_{0} / A_{0}} : \mcal{K}(B_{0} / L_{0}) \ar
\mcal{K}(A_{0} / L_{0})$
coincides with
$\sigma^{L_{0}}_{B_{0} / A_{0}} : \Omega^{n}_{B_{0} / L_{0}}
\ar \Omega^{n}_{A_{0} / L_{0}}$
of \cite{Ku} \S \textup{16}.
\end{lem}

\begin{proof}
Since $L_{0} \ar A_{0}$ is formally smooth, we get
$\mcal{K}(A_{0} / L_{0}) = \Omega^{n}_{A_{0} / L_{0}}$. The same for
$B_{0}$.
Consider the trivial DGA $L_{0}$. Then the universal
$B_{0}$-extension of $\Omega^{\bdot}_{A_{0} / L_{0}}$ is
$\Omega^{\bdot}_{B_{0} / L_{0}}$, so
$\sigma^{L_{0}}_{B_{0} / A_{0}}$ makes sense. To check that
$\sigma^{L_{0}}_{B_{0} / A_{0}} = \opn{Tr}_{B_{0} / A_{0}}$
we may reduce to the cases $A_{0} \ar B_{0}$ separable, or purely
inseparable of prime degree, and then use the properties of the trace.
\end{proof}

\begin{proof} (of the Theorem)\
Given any point $x \in X$ there is an open
neighborhood $U$ of $x$ which admits a factorization $\pi|_{U} = h g f$,
with $f : U \ar Y$ an open immersion; $g : Y \ar Z$ finite;
and $h : Z \ar S$ smooth of relative dimension $n$
(in fact one can take $Z$ open in $\mbf{A}^{n} \times S$).
This follows from quasi-normalization (\cite{Ku} Thm.\ B20)
and Zariski's Main Theorem (\cite{EGA} IV 8.12.3; cf.\
\cite{Ku} Thm.\ B16).
We can also assume $Y, Z, S$ are affine, say
$Y = \opn{Spec} B$,
$Z = \opn{Spec} A$ and $S = \opn{Spec} L$.
Let us write
$\tilde{\omega}^{n}_{B / L} := \Gamma(Y, \tilde{\omega}^{n}_{Y/S})$
and
$\mcal{K}^{\bdot}_{B / L} := \Gamma(Y, \mcal{K}^{\bdot}_{Y/S})$.
Also let us write
$B_{0} := k(Y)$, $A_{0} := k(Z)$ and $L_{0} := k(S)$.

By \cite{KW} \S 4,
\[ \tilde{\omega}^{n}_{B / L} =
\{ \beta \in \Omega^{n}_{B_{0} / L_{0}} \mid
\sigma^{L_{0}}_{B_{0} / A_{0}} (b \beta) \in \Omega^{n}_{A / L}
\text{ for all } b \in B \} . \]
One has
\[ \mcal{K}^{-n}_{B / L} = \mcal{K}(B_{0} / L_{0}) =
\Omega^{n}_{A_{0} / L_{0}}  \]
and the same for $A$.
According to Prop.\ \ref{prop6.4} there is a quasi-isomorphism
$\Omega^{n}_{A / L}[n]$ \linebreak
$ \ar \mcal{K}^{\bdot}_{A / L}$.
From the commutative diagram
\[ \begin{CD}
0 @> >> \mrm{H}^{-n} \mcal{K}^{\bdot}_{B / L} @> >>
\mcal{K}^{-n}_{B / L} @>{\delta}>>
\mcal{K}^{-n+1}_{B / L} \\
% @VVV @VVV @V{\opn{Tr}_{g}}VV @V{\opn{Tr}_{g}}VV \\
& & @VVV @V{\opn{Tr}_{g}}VV @V{\opn{Tr}_{g}}VV \\
0 @>>> \Omega^{n}_{A / L}  @>>>
\mcal{K}^{-n}_{A / L} @>{\delta}>>
\mcal{K}^{-n+1}_{A / L}
\end{CD} \]
and the isomorphism
\[ \mcal{K}^{-n+1}_{B / L} \cong
\opn{Hom}_{A} ( B, \mcal{K}^{-n+1}_{A / L}) \]
induced by $\opn{Tr}_{g}$
we conclude that
$\tilde{\omega}^{n}_{B / L} = \mrm{H}^{-n} \mcal{K}^{\bdot}_{B / L}$.
Since $\tilde{\omega}^{n}_{Y / S}$ and
$\mrm{H}^{-n} \mcal{K}^{\bdot}_{Y / S}$
are coherent sheaves and $f : U \ar Y$ is an open immersion, this shows
that
$\tilde{\omega}_{U/S} = \mrm{H}^{-n} \mcal{K}^{\bdot}_{U/S}$.
\end{proof}

\begin{cor}
If $X$ is a Cohen-Macaulay scheme then the sequence
\[ 0 \ar \tilde{\omega}^{n}_{X / S} \ar
\mcal{K}^{-n}_{X / S} \ar \cdots \ar
\mcal{K}^{m}_{X / S} \ar 0 \]
\textup{(}$m = \opn{dim} S$\textup{)} is exact.
\end{cor}

\begin{proof}
$X$ is Cohen-Macaulay iff any dualizing complex has a single nonzero
cohomology sheaf.
\end{proof}

\begin{exa} \label{exa7.1}
Suppose $X$ is an $(n+1)$-dimensional integral scheme and
$\pi : X \ar \opn{Spec} \mbb{Z}$ is a finite type
dominant morphism (i.e.\ $X$ has mixed characteristics). Then
$\pi$ is flat, equidimensional of dimension $n$ and
generically smooth. So
\[ \tilde{\omega}^{n}_{X / \mbb{Z}} =
\mrm{H}^{-n} \mcal{K}^{\bdot}_{X / \mbb{Z}} \subset
\Omega^{n}_{k(X) / \mbb{Q}} . \]
\end{exa}

\begin{rem} \label{rem7.7}
In the situation of Thm.\ \ref{thm7.4} there is a homomorphism
\[ \mrm{C}_{X} : \Omega^{n}_{X / S} \ar \mcal{K}^{-n}_{X / S} \]
called the {\em fundamental class of} $X / S$. According to \cite{KW},
when $\pi$ is flat one has
$\mrm{C}_{X}(\Omega^{n}_{X / S}) \subset \tilde{\omega}^{n}_{X/S}$; so
$\mrm{C}_{X} : \Omega^{n}_{X / S}[n] \ar \mcal{K}^{\bdot}_{X / S}$
is a homomorphism of complexes.
\end{rem}

\begin{rem}
In \cite{LS2} Theorem 11.2 we find a stronger statement
than our Theorem \ref{thm7.4}: $S$ is only
required to be an excellent equidimensional scheme without embedded
points, satisfying Serre's condition $\mrm{S}_{2}$; and $\pi$ is
finite type, equidimensional and generically smooth. Moreover,
for $\pi$ proper, the trace is compared to the integral
of \cite{HS} (cf.\ Remark \ref{rem8.1}).
The price of this generality is that the proofs in
\cite{LS2} are not self-contained but rely on rather complicated
results from other papers.
\end{rem}

% ** section 8 **

\section{The Isomorphism $\mcal{K}^{\bdot}_{X / S}
\protect \cong \pi^{!} \mcal{O}_{S}$}

In this section we describe the canonical isomorphism between the
complex $\mcal{K}^{\bdot}_{X / S}$ constructed in
Section 6, and the twisted inverse image $\pi^{!} \mcal{O}_{S}$
of \cite{RD}.
Recall that for residual complexes there is an inverse image
$\pi^{\triangle}$, and
$\pi^{\triangle} \mcal{K}^{\bdot}_{S / S} =
\mrm{E} \pi^{!} \mcal{O}_{S}$,
where $\mrm{E}$ is the Cousin functor corresponding to the dualizing
complex $\pi^{!} \mcal{O}_{S}$.
For an $S$-morphism $f : X \ar Y$ denote by
$\opn{Tr}^{\mrm{RD}}_{f}$ the homomorphism of graded sheaves
\[ \opn{Tr}^{\mrm{RD}}_{f} :
f_{*} \pi^{\triangle}_{X} \mcal{K}^{\bdot}_{S / S} \cong
f_{*} f^{\triangle} \pi^{\triangle}_{Y} \mcal{K}^{\bdot}_{S / S} \ar
\pi^{\triangle}_{Y} \mcal{K}^{\bdot}_{S / S} \]
of \cite{RD} Section VI.4.

\begin{thm} \label{thm8.10}
Let $\pi : X \ar S$ be a finite type morphism. Then there exists a
unique isomorphism of complexes
\[ \zeta_{X} : \mcal{K}^{\bdot}_{X / S} \ar
\pi^{\triangle} \mcal{K}^{\bdot}_{S / S} \]
such that for every morphism $f : X \ar Y$ the diagram
\begin{equation} \label{eqn8.4}
\begin{CD}
f_{*} \mcal{K}^{\bdot}_{X / S} @>{\opn{Tr}_{f}}>>
\mcal{K}^{\bdot}_{Y / S} \\
@V{f_{*} (\zeta_{X})}VV @V{\zeta_{Y}}VV \\
f_{*} \pi^{\triangle}_{X} \mcal{K}^{\bdot}_{S / S}
@>{\opn{Tr}^{\mrm{RD}}_{f}}>>
\pi^{\triangle}_{Y} \mcal{K}^{\bdot}_{S / S}
\end{CD}
\end{equation}
is commutative.
\end{thm}

The proof of Thm.\ \ref{thm8.10} is given later in this section,
after some preparation. Here is one corollary:

\begin{cor} \label{cor8.1}
If $f : X \ar Y$ is proper then $\opn{Tr}_{f}$ is a homomorphism of
complexes, and for any
$\mcal{M}^{\bdot} \in \msf{D}^{-}_{\mrm{qc}}(\msf{Mod}(X))$
the induced morphism
\[ f_{*} \mcal{H}om_{X}(\mcal{M}^{\bdot}, \mcal{K}^{\bdot}_{X / S}) \ar
\mcal{H}om_{X}(\mrm{R} f_{*} \mcal{M}^{\bdot}, \mcal{K}^{\bdot}_{Y / S}) \]
is an isomorphism.
\end{cor}

\begin{proof}
Use \cite{RD} Theorem VII.2.1 and Corollary VII.3.4.
\end{proof}

\begin{rem} \label{rem8.1}
In \cite{Hg3} and \cite{LS2} the authors prove that in their
respective constructions the trace
$\opn{Tr}_{f} : f_{*} f^{!} \mcal{N}^{\bdot} \ar \mcal{N}^{\bdot}$
is a homomorphism of complexes for any proper morphism $f$
and residual (resp.\ Cousin) complex $\mcal{N}^{\bdot}$
(cf.\ Remark \ref{rem6.6}).
\end{rem}

Let $Y = \opn{Spec} A$ be an affine noetherian scheme,
$X := \mbf{A}^{n} \times Y =$ \newline
$\opn{Spec} A \sqbr{t_{1}, \ldots, t_{n}}$
and $f : X \ar Y$ the projection.
Fix a point $x \in X$, and let $y := f(x)$,
$Z_{0} := \overline{\{x\}}_{\mrm{red}}$. Assume $Z_{0} \ar Y$ is finite.

\begin{lem} \label{lem8.10}
There exists an open set $U \subset Y$ containing $y$ and
a flat finite morphism $g : Y' \ar U$ s.t.:
\begin{enumerate}
\rmitem{i} $g^{-1}(y)$ is one point, say $y'$.
\rmitem{ii} Define $X' := \mbf{A}^{n} \times Y'$, and let
$f' : X' \ar Y'$, $h : X' \ar X$. Then for every point $x' \in h^{-1}(x)$
there is some section $\sigma_{x'} : Y' \ar X'$ of $f'$ with
$x' \in \sigma_{x'}(Y')$.
\end{enumerate}
\end{lem}

\begin{proof}
Choose any finite normal field extension $K$ of $k(y)$ containing $k(x)$.
Define recursively open sets
$U_{i} = \opn{Spec} A_{i} \subset Y$
and finite flat morphisms
$g_{i} : Y_{i} = \opn{Spec} A'_{i} \ar U_{i}$
s.t.\
$g_{i}^{-1}(y) = \{y_{i}\}$ and
$k(y_{i}) \subset K$,
as follows. Start with
$U_{0} = Y_{0} := Y$ and
$A'_{0} = A_{0} := A$.
If $k(y_{i}) \neq K$ take some $\bar{b} \in K - k(y_{i})$ and let
$\bar{p} \in k(y_{i}) \sqbr{t}$ be the monic irreducible polynomial of
$\bar{b}$. Choose a monic polynomial
$p \in \mcal{O}_{Y_{i}, y_{i}} \sqbr{t}$
lifting $\bar{p}$. There is some open set
$U_{i + 1} = \opn{Spec} A_{i + 1} \subset U_{i}$
s.t.\
$p \in (A'_{i} \otimes_{A_{i}} A_{i + 1}) \sqbr{t}$.
Define
$A'_{i + 1} := (A'_{i} \otimes_{A_{i}} A_{i + 1}) \sqbr{t} / (p)$
and $Y_{i + 1} = \opn{Spec} A'_{i + 1}$.
For $i = r$ this stops, and $k(y_{r}) = K$.

For every point
$x' \in \opn{Spec} (K \otimes_{k(y)} k(x))$
and $1 \leq i \leq n$ let $\bar{a}_{i, x'} \in k(x') \cong k(y_{r})$
be the image of $t_{i}$, and let
$a_{i, x'} \in \mcal{O}_{Y_{r}, y_{r}}$
be a lifting. Take an open set
$U = \opn{Spec} A_{r + 1} \subset U_{r}$
s.t.\ each
$a_{i, x'} \in A' = (A'_{r} \otimes_{A_{r}} A_{r + 1})$,
and define
$Y' := \opn{Spec} A'$.
So for each $x'$ the homomorphism
$B' = A' \sqbr{t} \ar A'$, $t_{i} \mapsto a_{i, x'}$ gives the desired
section $\sigma_{x'} : Y' \ar X'$.
\end{proof}

Let $Z_{i}$ be the $i$-th infinitesimal neighborhood of
$Z_{0}$ in $X$, so $f_{i} : Z_{i} \ar Y$ is a finite morphism.
Suppose we are given a quasi-coherent $\mcal{O}_{Y}$-module
$\mcal{M}$  which is supported on $\overline{\{ y \}}$. One has
\[ \mcal{H}^{n}_{Z_{0}} \left( \Omega^{n}_{X / Y} \otimes f^{*} \mcal{M}
\right) \cong
\lim_{i \ar} \mcal{E}xt^{n}_{X} \left( \mcal{O}_{Z_{i}},
\Omega^{n}_{X / Y} \otimes f^{*} \mcal{M} \right) \]
and by \cite{RD} Thm.\ VI.3.1
\[ \mcal{E}xt^{n}_{X} \left( \mcal{O}_{Z_{i}},
\Omega^{n}_{X / Y} \otimes f^{*} \mcal{M} \right) =
\mcal{H}^{0} f_{i}^{!} \mcal{M} . \]
Note that we can also factor $f_{i}$ through $\mbf{P}^{n} \times Y$,
so $f_{i}$ is projectively embeddable, and by \cite{RD} Thm.\ III.10.5
we have a map
\begin{equation} \label{eqn8.10}
\opn{Tr}_{f}^{\mrm{RD}} :
f_{*} \mcal{H}^{n}_{Z_{0}} \left( \Omega^{n}_{X / Y} \otimes
f^{*} \mcal{M} \right) \ar \mcal{M} .
\end{equation}

Now define
$\widehat{A} := \widehat{\mcal{O}}_{Y , y}$ and
$\widehat{B} := \widehat{\mcal{O}}_{X , x}$, with
$\mfrak{n} \subset \widehat{B}$ the maximal ideal and
$\phi = f^{*} : \widehat{A} \ar \widehat{B}$.
Set $M := \mcal{M}_{y}$, which is a discrete $\widehat{A}$-module.
We then have a natural isomorphism of $\widehat{A}$-modules
\begin{equation} \label{eqn8.1}
\left( f_{*} \mcal{H}^{n}_{Z_{0}} \left( \Omega^{n}_{X / Y} \otimes
f^{*} \mcal{M} \right) \right)_{y} \cong
\mrm{H}^{n}_{\mfrak{n}} \left(
\widehat{\Omega}^{n}_{\widehat{B} / \widehat{A}} \otimes_{\widehat{A}}
M \right) \cong \phi_{\#} M .
\end{equation}

\begin{lem} \label{lem8.1}
Under the isomorphism \textup{(\ref{eqn8.1})},
\[ \opn{Tr}_{f}^{\mrm{RD}} = \opn{Tr}_{\phi} : \phi_{\#} M \ar M . \]
\end{lem}

\begin{proof}
The proof is in two steps.\\
Step 1.\ Assume there is a section $\sigma : Y \ar X$ to $f$ with
$x \in W_{0} = \sigma(Y)$. The homomorphism
$\sigma^{*} : B = A \sqbr{\ul{t}} \ar A$ chooses
$a_{i} = \sigma^{*}(t_{i}) \in A$,
so after the linear change of variables $t_{i} \mapsto t_{i} - a_{i}$
we may assume that $\sigma$ is the $0$-section (i.e.\
$\mcal{O}_{W_{0}} = \mcal{O}_{X} / \mcal{O}_{X} \cdot \ul{t}$).
Let $W_{i}$ be the $i$-th infinitesimal neighborhood of $W_{0}$.
Since $f : W_{i} \ar Y$ is projectively embeddable, there is a trace
map
\[ \opn{Tr}_{f}^{\mrm{RD}} :
f_{*} \mcal{H}^{n}_{W_{0}} \Omega^{n}_{X / Y} \ar \mcal{O}_{Y} . \]
For any $a \in A$ one has
\begin{equation} \label{eqn8.11}
\opn{Tr}_{f}^{\mrm{RD}} \left(
\gfrac{a \mrm{d} t_{1} \wedge \cdots \wedge \mrm{d} t_{n}}
{t_{1}^{i_{1}} \cdots t_{n}^{i_{n}}} \right) =
\begin{cases}
a & \text{ if } \ul{i} = (1, \ldots, 1) \\
0 & \text{ otherwise} .
\end{cases}
\end{equation}
This follows from properties R6 (normalization) and
R7 (intersection) of the residue symbol (\cite{RD} Section III.9).
Alternatively this can be checked as follows. Note that
$\opn{Tr}_{f}^{\mrm{RD}}$ factors through
$\mrm{R} f_{*} \Omega^{n}_{\mbf{P}^{n}_{Y} / Y}$.
For the case $\ul{i} = (1, \ldots, 1)$ use \cite{RD}
Proposition III.10.1.
For $\ul{i} \neq (1, \ldots, 1)$ consider a change of coordinates
$t_{i} \mapsto \lambda_{i} t_{i}$, $\lambda_{i} \in A$. By
\cite{RD} Corollary III.10.2, $\opn{Tr}_{f}^{\mrm{RD}}$ is independent
of homogeneous coordinates, so it must be $0$.

Now since $W_{0} \cap f^{-1}(y) = Z_{0}$ we have
\[ \mcal{H}^{n}_{Z_{0}} \left( \Omega^{n}_{X / Y} \otimes f^{*} \mcal{M}
\right) \cong
\mcal{H}^{n}_{W_{0}} \left( \Omega^{n}_{X / Y} \otimes f^{*} \mcal{M}
\right)  \]
and so the formula for $\opn{Tr}_{f}^{\mrm{RD}}$ in (\ref{eqn8.10})
is given by (\ref{eqn8.11}).
But the same formula is used in \cite{Hg1} to define
$\opn{Tr}_{\phi}$.\\
Step 2.\
The general situation: take $g : Y' \ar Y$ as in Lemma \ref{lem8.10},
and set $Z_{0}' := Z_{0} \times_{Y} Y'$.
The flatness of $g$ implies there is a natural
isomorphism of $\mcal{O}_{Y'}$-modules
\[ g^{*} f_{*} \mcal{H}^{n}_{Z_{0}} \left( \Omega^{n}_{X / Y} \otimes
f^{*} \mcal{M} \right) \cong
f_{*}' \mcal{H}^{n}_{Z_{0}'} \left( \Omega^{n}_{X' / Y'} \otimes
{f'}^{*} \mcal{M}' \right) \]
(where $\mcal{M}' := g^{*} \mcal{M}$) and by \cite{RD} Thm.\ III.10.5
property TRA4 we have
\begin{equation} \label{eqn8.14}
g^{*}(\opn{Tr}_{f}^{\mrm{RD}}) = \opn{Tr}_{f'}^{\mrm{RD}} .
\end{equation}
Let $\widehat{A}' := \widehat{\mcal{O}}_{Y', y'} \cong A' \otimes_{A}
\widehat{A}$,
so $\widehat{A} \ar \widehat{A}'$ is finite flat. Therefore
\begin{equation} \label{eqn8.2}
 \widehat{A}' \otimes_{\widehat{A}} \mrm{H}^{n}_{\mfrak{n}}
\left( \widehat{\Omega}^{n}_{\widehat{B} / \widehat{A}}
\otimes_{\widehat{A}} M \right)
\cong
\bigoplus_{\mfrak{n}' \in Z_{0}'}
\mrm{H}^{n}_{\mfrak{n}'}
\left( \widehat{\Omega}^{n}_{\widehat{B}_{\mfrak{n}'} / \widehat{A}'}
\otimes_{\widehat{A}'} M' \right) .
\end{equation}
Here
$M' := \mcal{M}'_{y'} \cong \widehat{A}' \otimes_{\widehat{A}} M$
and
$\prod_{\mfrak{n}' \in Z_{0}'} \widehat{B}_{\mfrak{n}'}$
is the decomposition of $A' \otimes_{A} \widehat{B}$
to local rings. Write
$\phi'_{\mfrak{n}'} : \widehat{A}' \ar \widehat{B}_{\mfrak{n}'}$.
Direct verification shows that under the isomorphism (\ref{eqn8.2}),
\begin{equation} \label{eqn8.15}
1 \otimes \opn{Tr}_{\phi} = \sum_{\mfrak{n}' \in Z_{0}'}
\opn{Tr}_{\phi'_{\mfrak{n}'}} .
\end{equation}

Since $\widehat{A} \ar \widehat{A}'$ is faithfully flat it follows that
$M \ar M'$ is injective. In view of
the equalities (\ref{eqn8.14}) and (\ref{eqn8.15}), we conclude
that it suffices to check for each $\mfrak{n}' = x' \in Z_{0}$ that
$\opn{Tr}_{\phi'_{\mfrak{n}'}} =
\opn{Tr}^{\mrm{RD}}_{f'}$
on
$\mrm{H}^{n}_{\mfrak{n}'}
\left( \widehat{\Omega}^{n}_{\widehat{B}_{\mfrak{n}'} / \widehat{A}'}
\otimes_{\widehat{A}'} M' \right)$.
But there is a section $\sigma_{x'} : Y' \ar X'$, so we can apply
step 1.
\end{proof}

\begin{proof} (of Thm.\ \ref{thm8.10}.)\\
Step 1. (Uniqueness)\
Suppose
$\zeta_{X}' : \mcal{K}^{\bdot}_{X / S} \ar
\pi^{\triangle} \mcal{K}^{\bdot}_{S / S}$
is another isomorphism satisfying
$\opn{Tr}_{\pi} = \opn{Tr}^{\mrm{RD}}_{\pi} \pi_{*}(\zeta_{X}')$.
Then $\zeta_{X}' = a \zeta_{X}$ for some
$a \in \Gamma(X, \mcal{O}_{X}^{*})$, and by
assumption for any closed point $x \in X$ and
$\alpha \in \mcal{K}_{X / S}(x)$ there is equality
$\opn{Tr}_{\pi}(\alpha) = \opn{Tr}_{\pi} (a \alpha)$.
Now writing $s := \pi(x)$, it's known that
\[ \opn{Hom}_{\mcal{O}_{S, s}} \left(
\mcal{K}_{X / S}(x), \mcal{K}_{S / S}(s) \right) \]
is a free $\widehat{\mcal{O}}_{X, x}$-module with basis
$\opn{Tr}_{\pi}$. Therefore $a = 1$ in
$\widehat{\mcal{O}}_{X, x}$. Because this is true for all closed points
we see that $a = 1$.\\
Step 2.\
Assume $X = \mbf{A}^{n} \times S$ and $f = \pi$.
In this case there is a canonical isomorphism of complexes
\[ \mcal{K}_{X / S}^{\bdot} \cong \mrm{E} \Omega^{n}_{X / S}[n] \cong
\mrm{E} \pi^{!} \mcal{O}_{S} \cong
\pi^{\triangle} \mcal{K}_{S / S}^{\bdot}  \]
(cf.\ \cite{RD} Thm.\ VI.3.1 and our Prop.\ \ref{prop6.4}),
which we use to define
$\zeta_{X} : \mcal{K}^{\bdot}_{X / S} \ar
\pi^{\triangle} \mcal{K}^{\bdot}_{S / S}$.
Consider $x \in X$, $Z := \overline{\{x\}}_{\mrm{red}}$,
$s := \pi(x)$ and assume $x$ is closed in $\pi^{-1}(s)$.
By replacing $S$ with a suitable open neighborhood of $s$ we can
assume $Z \ar S$ is finite. Then we are allowed to apply Lemma
\ref{lem8.1} with $Y = S$, $\mcal{M} = \mcal{K}_{S / S}(s)$.
It follows that (\ref{eqn8.4}) commutes on
$\pi_{*} \mcal{K}_{X / S}(x) \subset \pi_{*} \mcal{K}^{\bdot}_{X / S}$.\\
Step 3.\ Let $X$ be any finite type $S$-scheme. For every affine open
subscheme $U \subset X$ we can find a closed immersion
$h : U \ar \mbf{A}^{n}_{S}$. Write $Y := \mbf{A}^{n}_{S}$
and let
$\pi_{U}$ and $\pi_{Y}$ be the structural morphisms.
Now
$\opn{Tr}_{h}$ induces an isomorphism
\[ \mcal{K}^{\bdot}_{U / S} \cong
\mcal{H}om_{Y}(\mcal{O}_{U}, \mcal{K}^{\bdot}_{Y / S}) , \]
and
$\opn{Tr}^{\mrm{RD}}_{h}$ induces an isomorphism
\[ \pi_{U}^{\triangle} \mcal{K}^{\bdot}_{S / S} \cong
\mcal{H}om_{Y}(\mcal{O}_{U},
\pi_{Y}^{\triangle} \mcal{K}^{\bdot}_{S / S}) . \]
So the isomorphism
$\zeta_{Y}$ of Step 2 induces an isomorphism
$\zeta_{U} : \mcal{K}^{\bdot}_{U / S} \ar
\pi_{U}^{\triangle} \mcal{K}^{\bdot}_{S / S}$, which satisfies
$\opn{Tr}_{\pi_{U}} = \opn{Tr}^{\mrm{RD}}_{\pi_{U}}
\pi_{U *}(\zeta_{U})$.
According to Step 1 the local isomorphisms $\zeta_{U}$ can be glued
to a global isomorphism $\zeta_{X}$.\\
Step 4.\ Let $f : X \ar Y$ be any $S$-morphism. To check (\ref{eqn8.4})
we may assume $X$ and $Y$ are affine, and in view of step 3
we may in fact assume
$Y = \mbf{A}^{m} \times S$ and
$X = \mbf{A}^{n} \times Y \cong \mbf{A}^{n + m} \times S$.
Now apply Lemma \ref{lem8.1} with $x \in X$ closed in its fiber
and
$\mcal{M} := \mcal{K}_{Y / S}(y)$.
\end{proof}

% ** references **

\end{document}